\documentclass[a4paper,11pt]{article}

\usepackage{graphicx}
\usepackage{epstopdf}
\usepackage{amssymb}
\usepackage{latexsym}

\usepackage{authblk}

\usepackage{geometry}
\usepackage{tensor}
\ifpdf
  \DeclareGraphicsExtensions{.eps,.pdf,.png,.jpg}
\else
  \DeclareGraphicsExtensions{.eps}
\fi
\usepackage{amsmath}
\usepackage{amssymb}

\usepackage{tikz}
\usepackage{tikz-cd}
\usepackage{pgfplots}

\DeclareMathOperator*{\outP}{\otimes}
\DeclareMathOperator*{\inP}{\cdot}
\newcommand{\oP}[1]{\outP\limits_{ #1}}
\newcommand{\iP}[1]{\inP\limits_{#1}}

\newcommand{\ie}{i.\,e.}%
\newcommand{\eg}{e.\,g.}%
\newcommand{\st}{s.\,t.}%
\newcommand{\wrt}{w.r.t.}%

\newcommand{\formComma}{\,\text{,}}
\newcommand{\formPeriod}{\,\text{.}}

\newcommand{\R}{\mathbb{R}}%

\newcommand{\ext}[1]{\underline{#1}}
\newcommand{\extTF}[1]{\underline{#1}}
\newcommand{\extNabla}{\underline{\nabla}\,}

\newcommand{\tensorField}[2]{\mathbb{T}^{#1}(#2)}

\newcommand{\landau}{\mathcal{O}}

\newcommand{\surf}{\mathcal{S}}

\newcommand{\manifold}{\mathcal{M}}
\newcommand{\normal}{\boldsymbol{\nu}}
\newcommand{\conormal}{\boldsymbol{n}}
\newcommand{\tangent}{\textup{T}}

\newcommand{\gbCmp}{g}
\newcommand{\gb}{\boldsymbol{\gbCmp}}%
\newcommand{\GbCmp}{G}%
\newcommand{\Gb}{\boldsymbol{\GbCmp}}%

\newcommand{\tb}{\boldsymbol{t}}%
\newcommand{\ub}{\boldsymbol{u}}%
\newcommand{\pb}{\boldsymbol{p}}%
\newcommand{\rb}{\boldsymbol{r}}%
\newcommand{\fb}{\boldsymbol{f}}%
\newcommand{\eb}{\boldsymbol{e}}%

\newcommand{\xb}{\boldsymbol{x}}%
\newcommand{\xbe}{\ext{\boldsymbol{x}}}%
\newcommand{\Xc}{\mathcal{X}}%
\newcommand{\qb}{\boldsymbol{q}}%
\newcommand{\psib}{\boldsymbol{\psi}}%

\newcommand*{\rom}[1]{\textup{\uppercase\expandafter{\romannumeral#1}}}

\newcommand{\shapeOperator}{\boldsymbol{B}}%
\newcommand{\gaussianCurvature}{\mathcal{K}}%
\newcommand{\meanCurvature}{\mathcal{H}}%
\newcommand{\Div}{\operatorname{div}}
\newcommand{\Rot}{\operatorname{Rot}}
\newcommand{\rot}{\operatorname{rot}}

\newcommand{\Christ}[2]{\Gamma^{#1}_{#2}}
\newcommand{\ChristE}[2]{\underline{\Gamma}^{#1}_{#2}}

\newcommand{\trace}{\operatorname{tr}_{\manifold}}%
\newcommand{\traceR}{\operatorname{tr}_3}%

\newcommand{\ProjSurf}{\operatorname{\Pi}}

\title{A finite element approach for vector- and tensor-valued surface PDEs}%

\author[1]{Michael Nestler\footnote{Corresponding author: michael.nestler@tu-dresden.de}}
\author[1]{Ingo Nitschke}
\author[1,2,3]{Axel Voigt}

\affil[1]{Institut f{\"u}r Wissenschaftliches Rechnen, Technische Universit{\"a}t Dresden, 01062 Dresden, Germany}
\affil[2]{Dresden Center for Computational Materials Science (DCMS), Technische Universit{\"a}t Dresden, 01062 Dresden, Germany}
\affil[3]{Center for Systems Biology Dresden (CSBD), Pfotenhauerstr. 108, 01307 Dresden, Germany}

\begin{document}

\maketitle

\begin{abstract}
We derive a Cartesian componentwise description of the covariant derivative of tangential tensor fields of any degree on general manifolds. This allows to reformulate any vector- and tensor-valued surface PDE in a form suitable to be solved by established tools for scalar-valued surface PDEs. We consider piecewise linear Lagrange surface finite elements on triangulated surfaces and validate the approach by a vector- and a tensor-valued surface Helmholtz problem on an ellipsoid. We experimentally show optimal (linear) order of convergence for these problems. The full functionality is demonstrated by solving a surface Landau-de Gennes problem on the Stanford bunny. All tools required to apply this approach to other vector- and tensor-valued surface PDEs are provided.
\end{abstract}

\section{Introduction}

Over the last decade a huge interest has developed for surface PDEs, in both, numerical analysis and applications. To deal with PDEs that are defined on curved surfaces requires to combine various mathematical disciplines, e.g., differential geometry, variational methods, numerics and even topology. But even if some numerical approaches exist, which rely on a geometric framework, mainly used in the context of computer graphics, see e.g. \cite{Craneetal_ACM_2013,DeGoesetal_ACM_2015}, the important breakthrough in the development of numerical methods for this type of PDEs is the avoidance of charts and atlases. Either the methods are based on a triangulated surface and require information on the surface solely through knowledge of the vertices, or an implicit surface representation is used and the problem is extended to the embedding space. Most of this work until today is concerned with scalar-valued surface PDEs, see \cite{Dziuk} for a review. In this case the coupling between the geometry of the surface and the PDE is weak and thus allows to solve these problems with small modifications of established numerical approaches in flat space. For vector- and tensor-valued surface PDEs this coupling becomes much stronger. Various new questions arise: How to define two tangential vectors on a triangulated surface to be parallel? What properties are required for a tangential surface tensor? What is an appropriate vector surface Laplacian? While the answers to these questions might be obvious in the context of differential geometry, the analysis of such highly nonlinear problems and the development of numerical methods are challenging. The research in this field is still at the very beginning, with some pioneering contributions derived for specific applications. In computer graphics tangent vector fields are essential in tasks such as texture synthesis and meshing and have been addresses by discrete differential geometry, see e.g. \cite{Hirani_ACM_2003,Craneetal_ACM_2013}. Only very few approaches exist, where these tools are used to solve physical problems, see e.g. \cite{Mohamedetal_JCP_2016,Nitschkeetal_book_2017,Nestleretal_JNS_2018} for a surface Navier-Stokes equation and a surface Frank-Oseen type model. Extensions to tensor fields of higher degree are rare and limited to specific tensor properties. In applications in geophysics surface PDEs on a sphere have a long tradition. These problems can be solved by spectral methods based on spherical harmonics expansions \cite{Backus_ARMA_1966,Barreraetal_EJP_1985,Freedenetal_MG_1995,Freedenetal_Springer_2009} and have recently also been considered for other applications, see e.g. \cite{Nestleretal_JNS_2018,Praetoriusetal_PRE_2018,Dunkeletal_PRL_2018} for a surface Frank-Oseen type model, a surface active phase-field-crystal model and active hydrodynamics, respectively. Also extensions to radial manifold shapes are possible, see e.g. \cite{Grossetal_JCP_2018}. But the approach is not suitable for general manifolds and also an extension to tensor fields of higher degree comes with various limitations. A numerical procedure applicable to general manifolds and a wide class of applications does not exist. However, there is a strong need for such a tool, resulting from increasing interests on vector- and tensor-valued surface PDEs. The interest results from the strong interplay between topology, geometry and dynamics. Examples are the location of topological defects in liquid crystals at extrema of the Gaussian curvature, the initiation of wrinkling in crystalline sheets at material defects or the possibility to control surface flow by changes in morphology. This interplay has the potential for breakthrough developments in materials science and synthetic biology, see \cite{Keberetal_Science_2014} for an example. Other applications are in general relativity. All ask for a general numerical approach for vector- and tensor-valued surface PDEs. We here provide an approach, which is based on a reformulation of the problems in Cartesian coordinates, which allows for a componentwise solution using tools for scalar-valued surface PDEs. Similar approaches for vector-valued surface PDEs have been considered in \cite{Nestleretal_JNS_2018,Hansboetal_arXiv_2017,Grossetal_arXiv_2017,Jahnkuhnetal_arXiv_2017,Olshanskiietal_arXiv_2018} for a surface Frank-Oseen type problem, a surface Laplace equation and a surface Stokes problem, respectively. However, all these approaches make strong assumptions on the extension of the surface quantities in the embedding space.
Here, we introduce an approach, which does not require these assumptions and which is applicable for vector- and tensor-valued surface PDEs. 
Starting with the manifold and the covariant derivative $\nabla$, we extend the manifold and the quantities on it to an embedding thin film with the corresponding embedding space/thin film derivative $\extNabla$. 
By using thin film coordinates we can separate tangential and normal contributions, which enables to identify the tangential parts with $\nabla$ and to express $\extNabla$ in terms of partial derivatives along the Euclidian basis. This reformulation maintains physical invariance and allows to rewrite any manifold bound, tensor-valued PDE into a set of scalar-valued PDE's along an Euclidian basis. In this form the problem can be solved for each component by established finite element methods. We numerically demonstrate optimal order of convergence for two examples and show illustrating results for 2-tensors on complex geometries. We further provide all necessary tools to solve general vector- and tensor-valued surface PDEs within a numerical framework for scalar-valued surface PDEs. 
The paper is structured as follows: In section \ref{sec:2} we provide the basic notations and geometric preliminaries, and derive a description of the covariant derivative in Cartesian coordinates for tensor fields of degree $d$. Based on this result a generic componentwise finite element formulation is presented in section \ref{sec:3}, which is demonstrated and validated on vector- and tensor-valued surface Helmholtz problems. Readers not interested in details of the derivation can skip the initial part of section \ref{sec:3} and focus on the discussion of the Cartesian description ot the covariant derivative starting at \eqref{eq:thinFilmCoordsNabla}. Convergence results, demonstrating the optimal order, are provided and the sensitivity of the error on various approximations is studied. To show the full functionality of the approach we solve a surface Landau-de Gennes model on the Stanford bunny in section \ref{sec:4}. Conclusions are drawn in section \ref{sec:5} and details of the derivation and extensions which allow to use the approach for general vector- and tensor-valued surface PDEs  are provided in the Appendices. 

\section{Notation and geometric preliminaries}
\label{sec:2}

\subsection{Tensor notation} An algebraic approach of modeling tensors is to consider them as outer products of vectors. The resulting tensor spaces conform to the concept of an algebraic vector space, enabling a description along basis and coordinates, e.g. $\tb = t^{ij} \eb_i \otimes \eb_j$ describing a 2-tensor $\tb$ by contravariant coordinates $t^{ij}$ and a covariant basis $\eb_i \otimes \eb_j$ formed by the outer product of the (local) basis vectors, see \cite{abraham2012manifolds,Schouten1954,Jaenich2013,kuhnel2006differential} for a general introduction.

We use Einstein sum convention for all coordinate descriptions. Further to distinguish between local/manifold bound coordinates and Cartesian/flat space coordinates we use two sets of indices of lower case and upper case letters, respectively. We focus on a description by covariant coordinates and are concerned with full contractions of tensor fields which we denote by ``:". Under this conditions covariant and contravariant descriptions are equivalent. 

Beside index notation we also use operator notation. To clarify the association of inner and outer products we also use subscripts to indicate the involved indices of the left side tensor. For the inner product ``$\inP$", e.g., we describe a contraction of the 3-tensor $\ub$ and the 2-tensor $\qb$ along the second components by
\begin{align}
 [\ub \iP{2} \qb]_{i_1 i_2 i_3} = \tensor{u}{_{i_1 }^{k}_{i_2}}\tensor{q}{_{i_3}_{k}}.
\end{align}
For the outer product ``$\outP$" we add as subscripts the indices of the resulting tensor, indicating where the left side tensor components are added, e.g., for the product of two 2-tensors $\qb$ and $\rb$ 
\begin{align}
 [\qb \oP{1,4}\rb]_{i_1 \hdots i_4} = \tensor{q}{_{i_1 i_4}} \tensor{r}{_{i_2 i_3}}.
\end{align}
In cases where the use of ``$\inP$" and ``$\outP$" is unambiguous we drop the subscripts for better readability.

\subsection{Local coordinates and covariant derivatives} We consider $n-1$ dimensional oriented Riemannian manifolds $\manifold$ embedded in $\R^n$, which are described along the chart-atlas by a set of local coordinates $\{u_i\} \subset \R^{n-1}$ such that $\Xc (\{u_i\}) \mapsto \xb \in \R^n$. The local basis vectors are given by the partial derivatives of $\Xc$: $\eb_i = \partial_i \Xc \in \tangent^1\manifold \subset \R^n$. Using a scalar product given by the full contraction of the tensor we define the Riemannian manifolds metric $\gb$
\begin{align}
\tensor{\gbCmp}{_{ij}} = \partial_i \Xc : \partial_j \Xc \;\text{, where } \; \tensor{\gbCmp}{^{ji}}=\tensor{[\gb^{-1}]}{^{ij}} ,\quad \tensor{\gbCmp}{_{i}^{j}} = \tensor{\delta}{_{i}^{j}}.
\end{align}
Let $\normal$ be the manifold normal. We define the shape operator by $\shapeOperator$, with $B_{ij} = -\partial_i \Xc : \partial_j \normal $.  This and its invariants $ \gaussianCurvature = \det\{\tensor{B}{^i_j}\} $ (Gaussian curvature) and $\meanCurvature=\tensor{B}{^i_i} $ (mean curvature)
can also be expressed in Cartesian coordinates of the embedding space by $B_{IJ} = - \tensor{\ProjSurf}{^{L}_{I}} \partial_L \nu_J$, with the surface identity $\tensor{\ProjSurf}{^{I}_{J}} = \tensor{\delta}{^{I}_{J}} - \nu^I\nu_J$ (in the sense of $\ProjSurf_{ij} = \gbCmp_{ij}$ and $\ProjSurf \inP \normal = \normal \inP \ProjSurf = 0 $). 
Note that we use lower case Latin indices for local manifold coordinates and capital Latin indices for local embedding space coordinates. On a manifold we consider covariante derivatives $\nabla$ which are metric compatible, namely $\nabla \gb = 0$ (the Levi-Civita connection). These can be described using the Christoffel symbols 
\begin{align}
\Christ{k}{ij} = \frac{1}{2}\tensor{\gbCmp}{^{kl}}\left( \partial_i \tensor{\gbCmp}{_{jl}} + \partial_j\tensor{\gbCmp}{_{il}} - \partial_l \tensor{\gbCmp}{_{ij}} \right)
\end{align}
and partial derivatives $\partial_i$. Depending on the tensorial degree $d $ of the covariant tensor field $\tb$ the covariant derivative, 
see \eg\ \cite{kuhnel2006differential, Schouten1954}, reads
\begin{align} \label{eq:covariant gradient}
[\nabla\tb]_{i_1 \hdots i_d k} = \;& \partial_k t_{i_1 \hdots i_d} - [ \Christ{l}{k i_1} t_{l i_2 \hdots i_d} +  \hdots + \Christ{l}{k i_d} t_{i_1 \hdots i_{d-1}l} ]\formPeriod
\end{align}

\subsubsection{Tensor valued fields on manifolds} Let $\tensorField{d}{\manifold} = \{ \tb \; :\; \manifold \mapsto \tangent^d\manifold \}$ denote the functionspace of smooth tangent tensor fields of degree $d$ on $\manifold$, with $\tangent^d \manifold = \cup_{\xb \in \manifold} \tangent^d_{\xb} \manifold$ and $\tangent^d_{\xb} \manifold$ the tangent space of d-tensors on $\manifold$ at $\xb \in \manifold$. We define an associated inner product by 
\begin{align}
\langle \tb, \qb  \rangle_d & = \int_{\manifold} \left( \tb(\xb), \qb(\xb) \right)_{d,\manifold} \mathrm{d}\manifold .
\end{align}
Where $\left( .,.\right)_{d,\manifold}$ denotes the canonical pointwise inner product, the full contraction of the tensors with respect to the manifold metric, e.g. for two 2-tensors $\left( \tb, \qb \right)_{2,\manifold} =  \tensor{t}{_{ij}}\tensor{q}{^{ij}}$. Please note, this coincides only for $\gbCmp_{ij} = \delta_{ij}$ with the previously introduced contraction ``:". Along the derivations we will also use non-tangential tensor fields $\ext{\tb}$ extended to an embedding domain. The projection to $\tangent^d\manifold$ is defined by
\begin{align}
\ProjSurf[\ext{\tb}]_{I_1 \hdots I_d}
    &= \tensor{\ProjSurf}{^{L_1 \hdots L_d}_{I_1 \hdots I_d}} \tensor{\ext{t}}{_{L_1 \hdots L_d}}
    = \tensor{\ProjSurf}{^{L_1}_{I_1}} \hdots \tensor{\ProjSurf}{^{L_d}_{I_d}} \tensor{\ext{t}}{_{L_1 \hdots L_d}}
\end{align}
and we extend the previous inner product modulo tangential parts in the arguments by
\begin{align} \label{eq:innerTan}
\langle \ext{\tb}, \ext{\qb}  \rangle_d & =  \langle \ProjSurf[\ext{\tb}], \ProjSurf[\ext{\qb}]  \rangle_d.
\end{align}

\subsection{Thin film extension}\label{sec:ThinFilmExtension} As an instructive link between the description in local manifold basis and embedding space basis we consider a sufficient smooth thin film continuation of surface quantities in a tubular extension of $\manifold$. For sufficient thin extensions an unique description by thin film coordinates exist which is given by augmenting the local basis by the manifolds oriented normal $\normal$ 
\begin{align}
\extTF{\xb} = \extTF{\mathcal{X}}(\{u_i\},\xi) = \Xc(\{u_i\}) + \xi\normal(\{u_i \}).
\end{align}
We extend all manifold bound quantities smoothly to the thin film, solemnly requiring $\ProjSurf[\extTF{\tb}|_{\xi=0}] = \tb$. The thin film metric is given by $\GbCmp_{IJ}=\partial_I\extTF{\mathcal{X}} : \partial_J\extTF{\mathcal{X}} $ and can be related to the manifold metric $\gb$ by
\begin{align}
& \mbox{symbolic:} && \mbox{thin film coordinates:} \nonumber \\
&\ProjSurf[\Gb] = \gb - 2\xi\shapeOperator + \xi^2\shapeOperator^2 && \tensor{G}{_{ij}} = \tensor{g}{_{ij}} - 2\xi \tensor{B}{_{ij}} + \xi^2 \tensor{B}{_{i}^{k}}\tensor{B}{_{kj}} \\
&\ProjSurf \iP{1} \Gb \iP{2}  \normal = \normal \iP{1} \Gb \iP{2}  \ProjSurf = 0 &&  \tensor{G}{_{i\xi}} = \tensor{G}{_{\xi i}} = 0 \nonumber \\
&\normal \iP{1} \Gb \iP{2}  \normal = 1 && \tensor{G}{_{\xi\xi}} = 1 \nonumber
\end{align}
The thin film Christoffel symbols $\ChristE{}{}$ and the thin film covariant derivative $\extNabla$ can be expressed in thin film coordinates,
\ie\ for indices  $\{I,J, \hdots \} = \{i,j, \hdots ,\xi\}$, and read
\begin{align} 
    \ChristE{K}{IJ} &= \frac{1}{2}\tensor{\GbCmp}{^{KL}}\left( \partial_I \tensor{\GbCmp}{_{JL}} + \partial_J\tensor{\GbCmp}{_{IL}} - \partial_L \tensor{\GbCmp}{_{IJ}} \right)\formComma \notag\\
    [\extNabla\extTF{\tb}]_{I_1 \hdots I_d K} &= \partial_K \extTF{t}_{I_1 \hdots I_d} - [ \ChristE{L}{K I_1} \extTF{t}_{L I_2 \hdots I_d} +  \hdots + \ChristE{L}{K I_d} \extTF{t}_{I_1 \hdots I_{d-1}L} ] \label{eq:ex_covariant gradient}
\end{align}
We consider a Taylor expansion from the surface in normal direction to relate $\ChristE{}{}$ to $\Christ{}{}$, for details see \ref{sec:AppendixThinFilmChristoffelSymbols},
\begin{align}
\ChristE{k}{ij} = \Christ{k}{ij} + \landau(\xi), \quad \ChristE{\xi}{ij} = \tensor{B}{_{ij}} + \landau(\xi), \quad \ChristE{k}{i\xi} = \ChristE{k}{\xi i} = -\tensor{B}{_{i}^{k}} + \landau(\xi), \quad
\ChristE{K}{\xi\xi} = \ChristE{\xi}{K\xi} = \ChristE{\xi}{\xi K} = 0. 
\end{align}

\subsection{Cartesian description for covariant derivatives of tensor fields on manifolds} 
The Cartesian description of eq. (\ref{eq:covariant gradient}) is developed in two parts.
Firstly, we relate $ \extNabla $ in \eqref{eq:ex_covariant gradient} with $\nabla$ in \eqref{eq:covariant gradient} regarding the conformities of 
local thin film and manifold coordinate systems.
Secondly, we use the principle of covariance for $\extNabla $ and tensor calculus in the embedding space. This e.g. allows to express the covariant derivative \wrt\ Cartesian coordinates.

We therefore consider thin film extensions of the tangential fields $\tb$, $\nabla\tb$, $\shapeOperator$ and non-tangential fields $\ext{\tb}$, $\extNabla\ext{\tb}$, $\normal$.
The tangential part of the flat space gradient can be described by thin film coordinates and its tangential directions $i_1, \hdots, i_d, k$ as
\begin{align}
[\extNabla\ext{\tb}]_{i_1 \hdots i_d k} = \partial_k \ext{t}_{i_1 \hdots i_d} - \left[ \ChristE{L}{k i_1} \ext{t}_{L i_2 \hdots i_d} +  \hdots + \ChristE{L}{k i_d} \ext{t}_{i_1 \hdots i_{d-1}L} \right].
\end{align}
Separating the summation in the Christoffel symbols in tangential $i_1 \hdots i_d$ and normal direction $\xi$ leads to
\begin{align}
[\extNabla\ext{\tb}]_{i_1 \hdots i_d k} =  \partial_k \ext{t}_{i_1 \hdots i_d}  - [ \ChristE{l}{k i_1} \ext{t}_{l i_2 \hdots i_d}  +  \hdots + \ChristE{l}{k i_d} \ext{t}_{i_1 \hdots i_{d-1}l}  ] - [ \ChristE{\xi}{k i_1} \ext{t}_{\xi i_2 \hdots i_d} +  \hdots + \ChristE{\xi}{k i_d} \ext{t}_{i_1 \hdots i_{d-1}\xi} ].
\end{align}
On the other hand, the extended covariant gradient $\nabla \tb$ expressed in thin film coordinates reads
\begin{align}
[\nabla\tb]_{i_1 \hdots i_d k} =  \partial_k t_{i_1 \hdots i_d}  - [ \ChristE{l}{k i_1} t_{l i_2 \hdots i_d} +  \hdots + \ChristE{l}{k i_d} t_{i_1 \hdots i_{d-1}l} ]  - [ \ChristE{\xi}{k i_1} \underbrace{t_{\xi i_2 \hdots i_d}}_{=0} +  \hdots + \ChristE{\xi}{k i_d} \underbrace{t_{i_1 \hdots i_{d-1}\xi}}_{=0} ].
\end{align}
At the manifold ($\xi=0$) and due to $\ProjSurf[\extTF{\tb}|_{\xi=0}] = \tb$ both expressions can be related to each other and we obtain 
\begin{align}
[\extNabla\ext{\tb}]_{i_1 \hdots i_d k} = [\nabla\tb]_{i_1 \hdots i_d k} - [ \ChristE{\xi}{k i_1} \ext{t}_{\xi i_2 \hdots i_d} +  \hdots + \ChristE{\xi}{k i_d} \ext{t}_{i_1 \hdots i_{d-1}\xi} ].
\end{align}
Identifying the thin film Christoffel symbols in normal direction by the shape operator, $\ChristE{\xi}{k i_d} = \tensor{B}{_{k i_d}}$, we obtain 
\begin{align}
[\nabla\tb]_{i_1 \hdots i_d k} = [\extNabla\ext{\tb}]_{i_1 \hdots i_d k} + [\tensor{B}{_{k i_1}} \ext{t}_{\xi i_2 \hdots i_d} +  \hdots + \tensor{B}{_{k i_d}} \ext{t}_{i_1 \hdots i_{d-1}\xi} ],
\end{align}
which can also be written in Cartesian coordinates,
\begin{align}
\label{eq:thinFilmCoordsNabla}
[\nabla\tb]_{I_1 \hdots I_d K} &=  \ProjSurf[\extNabla\ext{\tb}]_{I_1 \hdots I_d K} 
    + [\tensor{B}{_{K I_1}} \tensor{\ProjSurf}{^{J_2 \hdots J_d}_{I_2 \hdots I_d}} \tensor{\ext{t}}{^{L}_{J_2 \hdots J_d}} 
        +  \hdots + \tensor{B}{_{K I_d}}\tensor{\ProjSurf}{^{J_1 \hdots J_{d-1}}_{I_1 \hdots I_{d-1}}} \tensor{\ext{t}}{_{J_1 \hdots J_{d-1}}^{L}} ] \nu^{L}\formComma
\end{align}
since lower case thin film indices refer to the tangential parts of a tensor.
Reviewing this result we note that the Christoffel symbols, and therefore any explicit dependence on the surface metric $\gb$ and its smooth description, has vanished. Only the quantities $\normal$ and $\shapeOperator$ remain, making the description independent from the assumed chart-atlas mechanism. Also the smooth thin film extension is not used explicitly anymore, solemnly its existence is required. Further note that this description is exact, therefore retaining the tangential property of the covariant gradient, namely
\begin{align}
\nabla\tb = \ProjSurf[\extNabla\ext{\tb}] + \sum_{m=1}^{d} \shapeOperator \oP{m,d+1} \ProjSurf\left[\ext{\tb} \iP{m} \normal \right]  \; \in \; \tensorField{d+1}{\manifold} \quad \forall \; \ext{\tb}\mbox{, such that } \ProjSurf[\ext{\tb}] = \tb.
\label{eqn:exactCovarinatDerivative}
\end{align}
We see that $\nabla\tb = \ProjSurf[\extNabla (\ProjSurf [\ext{\tb}])] $ is valid by applying the product rule several times. However, we want to highlight here, that this simple and natural seeming result is not a definition, as in \cite{Hansboetal_arXiv_2017}, but rather the consequence of the established covariant derivative \eqref{eq:covariant gradient}, which does not need an embedding space, and the considered embedding of $\manifold$ and extention of $\tensorField{d}{\manifold}$. In particular we obtain for 0-tensors (scalars), 1-tensors (vectors), 2-tensors, $\ldots$
\begin{align}
d=0,&\, s:\manifold \mapsto  \tangent^0(\manifold)   & [\nabla s]_i &= \partial_i s                                                                               & [\nabla s] &= \Pi[\extNabla \ext{s}] \nonumber \\
d=1,&\, \pb:\manifold \mapsto \tangent^1(\manifold) & [\nabla \pb]_{ij} &= \partial_j p_i  - \Christ{k}{ij} p_k                                        & [\nabla \pb] &=  \Pi[\extNabla \ext{\pb}] + (\ext{\pb} \iP{} \normal) \shapeOperator  \nonumber\\
d=2,&\, \qb:\manifold \mapsto \tangent^2(\manifold) & [\nabla \qb]_{ijk} &= \partial_k q_{ij}  - \Christ{l}{ki} q_{lj} - \Christ{l}{kj} q_{il} & [\nabla \qb] &=  \Pi[\extNabla \ext{\qb}]  + \shapeOperator \oP{1,3} (\normal \iP{} \ext{\qb} \iP{2}\ProjSurf ) 
                          + \shapeOperator \oP{2,3} (\ProjSurf \iP{} \ext{\qb} \iP{2} \normal) \nonumber\\
\vdots & & & \nonumber
\end{align}
While the description of the covariant derivative for scalar quantities as $[\nabla s] = \Pi[\extNabla \ext{s}] = \extNabla \extTF{s} - \normal(\normal\cdot \extNabla \extTF{s})$ is well established, see e.g. Definition 2.3 in  \cite{Dziuk}, the description in Cartesian coordinates for tensors depends on the degree $d$ and requires additional coupling terms with geometric quantities, which increase with the degree $d$. 

The reformulations for the other first-order operators $\Div \tb$ and $\Rot \tb$ are derived along the manifold metric $\gb$ and the compatible Levi-Civita tensor $\boldsymbol{E}$ and are derived in \ref{sec:AppendixDifferentialOperators}.

With these expressions we have all tools available to rewrite general vector- and tensor-valued surface PDEs in Cartesian coordinates and solve them for each component by established methods for scalar-valued quantities, such as the surface finite element method \cite{DziukElliott_JCM_2007,DziukElliott_IMAJNA_2007,Dziuk}, level set approaches \cite{Bertalmioetal_JCP_2001,Greeretal_JCP_2006,Stoeckeretal_JIS_2008,Dziuketal_IFB_2008}, diffuse interface approximations \cite{Raetzetal_CMS_2006}, trace- and cut-finite element methods \cite{Olshanskiietal_LNCSE_2018,Burmanetal_IJNME_2015} or any other method. 

\section{Generic componentwise finite element formulation}
\label{sec:3}

The Cartesian descriptions can now be used to derive a generic componentwise finite element formulation for vector- and tensor-valued surface PDEs. For 1-tensors (vectors) similar approaches have been considered for specific applications. However, beside \cite{Hansboetal_arXiv_2017}, all use $[\nabla \pb] = \Pi[\extNabla \ext{\pb}]$, which only holds if $\normal \cdot \ext{\pb} = 0$. In \cite{Nestleretal_JNS_2018,Reutheretal_PF_2018} this is enforced by a penalty approach and in \cite{Jahnkuhnetal_arXiv_2017} by using a Lagrange multiplier. We will compare these approaches and analyze the relevance of the additional terms in \eqref{eqn:exactCovarinatDerivative} at the end of this section. 

\subsection{Variational formulation for vector- and tensor-valued surface Helmholtz problems}
To demonstrate the approach we consider a Helmholtz problem on a manifold without boundary embedded in $\R^3$. 
\begin{align}
-\Div \left(  \nabla \tb \right) + \tb = \fb \quad \mbox{on }\manifold
\label{eqn:basicLaplacePDE}
\end{align}
For the treatment of boundary terms we refer to \ref{sec:AppendixComponentwiseBoundaryConditions}. The variational formulation reads: Find $\tb \in {\mathbb{H}}^{d,1}_{tan}(\manifold)$, s.t. 
\begin{align}\label{eq:weakform}
a(\tb,\psib) = l(\psib) \quad \forall \psib \in {\mathbb{H}}^{d,1}_{tan}(\manifold)
\end{align}
with 
\begin{align}
a(\tb,\psib) = \langle \nabla \tb, \nabla \psib \rangle_{d+1} + \langle \tb, \psib  \rangle_d, \qquad
l(\psib) = \langle \fb, \psib  \rangle_d
\end{align}
with ${\mathbb{H}}^{d,1}_{tan}(\manifold)$ the Sobolev space of tangent tensor fields of degree $d$ on $\manifold$ with covariant gradient and an appropriate norm. Using the derived identities, this can also be written in the extended space and the problem expressed in Cartesian coordinates.
This extended formulation reads: Find $\ext{\tb} \in [H^1(\manifold)]^{n,d}$, s.t.
\begin{align}
\label{eq:extweakform}
a(\ext{\tb},\ext{\psib}) = l(\ext{\psib}) \quad \forall \ext{\psib} \in [H^1(\manifold)]^{n,d}
\end{align}
with 
\begin{align}\label{eq:extbiandlinform}
a(\ext{\tb},\ext{\psib}) = \langle \extNabla\ext{\tb} + \sum_{m=1}^{d} \shapeOperator \oP{d+1,m} \left(\ext{\tb} \iP{m} \normal \right), \extNabla\ext{\psib} + \sum_{m=1}^{d} \shapeOperator \oP{d+1,m} \left(\ext{\psib} \iP{m} \normal \right) \rangle_{d+1} + \langle \ext{\tb}, \ext{\psib}  \rangle_d, \qquad
l(\psib) = \langle \ext{\fb}, \ext{\psib}  \rangle_d
\end{align}
and $[H^1(\manifold)]^{n,d}$ the embedded product Sobolev space of tensor fields on $\manifold$ defined by
\begin{align}
[H^1(\manifold)]^{n,d} = \{ \ext{\tb} \in [L^2(\manifold)]^{n,d} \; | \; \extNabla\ext{\tb} \iP{d+1 } \ProjSurf \in [L^2(\manifold)]^{n,d+1}  \}, \quad \ext{\tb} = \mathbb{M}_{I_1 \hdots I_d} \eb^{I_1}\outP \hdots \outP \eb^{I_d} \; \forall \ext{\tb} \in [H^1(\manifold)]^{n,d}
\end{align}
with Cartesian basis $\eb^{I_1}\outP \hdots \outP \eb^{I_d}$ and the coefficient function matrix $\mathbb{M}_{I_1 \hdots I_d}(\xbe)$. 
The space $ [H^1(\manifold)]^{n,d} $ contains the variational space of \eqref{eq:extweakform} such that $ [H^1(\manifold)]^{n,d}  \supset {\mathbb{H}}^{d,1}_{tan}(\manifold) \cong \ProjSurf\left[ [H^1(\manifold)]^{n,d} \right] $ and is orthogonal decomposable into 
$ [H^1(\manifold)]^{n,d} = {\mathbb{H}}^{d,1}_{tan}(\manifold) \oplus {\mathbb{H}}^{d,1}_{nor}(\manifold) $ \wrt\ projection operator $ \ProjSurf $,
\ie\ for all $ \ext{\tb} = \tb + \tb_{nor} \in [H^1(\manifold)]^{n,d} $ holds
$ \tb = \ProjSurf[\ext{\tb}] \bot \tb_{nor} = \ext{\tb} -\ProjSurf[\ext{\tb}]  $, $ \tb\in{\mathbb{H}}^{d,1}_{tan}(\manifold) $ and
$ \tb_{nor}\in{\mathbb{H}}^{d,1}_{nor}(\manifold)   $.
Since $ {\mathbb{H}}^{d,1}_{nor}(\manifold) $ is located in the kernel of the symmetric bilinear form $ a $ in 
\eqref{eq:extbiandlinform}, the solution of \eqref{eq:extweakform} is not unique.
To overcome this issue and preserve also the well-posedness of problem  \eqref{eq:weakform} in the bigger space 
$  [H^1(\manifold)]^{n,d} \supset {\mathbb{H}}^{d,1}_{tan}(\manifold) $,  
 the problem has to be modified:  
Find $\ext{\tb} \in [H^1(\manifold)]^{n,d}$, s.t. 
\begin{align}
\label{eq:modextweakform}
\tilde{a}(\ext{\tb},\ext{\psib}) = l(\ext{\psib}) \quad \forall \ext{\psib} \in [H^1(\manifold)]^{n,d}
\end{align}
with 
\begin{align}
\tilde{a}(\ext{\tb},\ext{\psib}) = a(\ext{\tb},\ext{\psib}) + \langle p(\ext{\tb}), \ext{\psib}  \rangle
\end{align}
where $p$ is an appropriate penalization function s.t. 
$ \langle p(\ext{\tb}), \ext{\psib} \rangle = c \langle \ext{\tb} - \ProjSurf[\ext{\tb}] , \ext{\psib} \rangle  $ for 
$ c\in\R\!\setminus\!\{0\} $
and all 
$ \ext{\tb},\ext{\psib} \in  [H^1(\manifold)]^{n,d} $.
As a consequence, problem \eqref{eq:modextweakform} is equivalent to \eqref{eq:weakform} and the additional problem:
Find $ \tb_{nor}\in{\mathbb{H}}^{d,1}_{nor}(\manifold) $, \st
\begin{align}
    \langle \tb_{nor}, \psib_{nor} \rangle = 0 \quad \forall \psib_{nor}\in{\mathbb{H}}^{d,1}_{nor}(\manifold)
\end{align}
where $ \ext{\tb} = \tb + \tb_{nor} $.
Both problems are independent and well-posed and therefore also \eqref{eq:modextweakform}.

To validate this approach we will consider two cases: a 1-tensor (vector) formulation, which is related to the Maxwell equations in electromagnetic theory and a 2-tensor formulation, where we enforce the tensor field to be symmetric and traceless, which leads to a problem related to the Einstein equations in general relativity. However, we are here not interested in these applications, but only in the structure of the equation, which allows to construct an analytical solution by considering a smooth tensor field $\tb^*$ and defining $\fb = -\Div \left(  \nabla \tb^* \right) + \tb^*$. We will construct appropriate surface finite element discretizations for these 
problems based on \eqref{eq:modextweakform}. To account for the special tensor properties in the 2-tensor formulation will require some further modifications.

\subsection{Vector-valued surface Helmholtz problem}
\begin{figure}[ht]
\begin{center}
\includegraphics[height=0.3\textwidth]{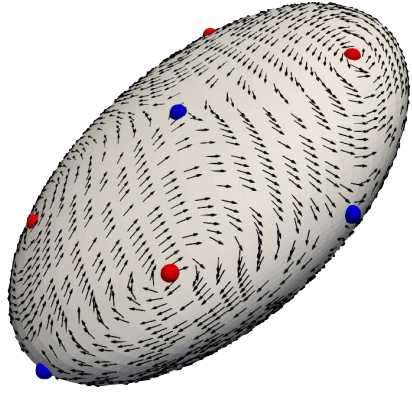}
\end{center}
\caption{\textbf{Normalized analytical solution $\pb^*$ with 8 vorticies (topological charge $+1$, red) and 6 saddle points (topological charge $-1$, blue). (colors online)}}
\label{fig:vectorHelmholtz}
\end{figure}

We consider a surface vector-valued Helmholtz problem, which reads in weak form 
\begin{align}
\tilde{a}(\ext{\pb},\ext{\psib}) = l(\ext{\psib}) \quad \forall \ext{\psib} \in [H^1(\manifold)]^{3,1}
\end{align}
with $\manifold$ to be an ellipsoid with major axis $A=1$, $B=0.5$ and $C=1.5$ and $\fb = -\Div \left(  \nabla \pb^* \right) + \pb^*$, with $\pb^* = \Rot(xyz)$, see figure  \ref{fig:vectorHelmholtz}. 
The solution contains various defects with topological charge (winding number) $\pm 1$, while the sum of these charges (winding numbers) is two, equal to the Euler characteristic of the ellipsoid, thus fulfilling the Poincare-Hopf theorem. 

\begin{figure}[ht]
\begin{center}
\includegraphics[height=0.2\textwidth]{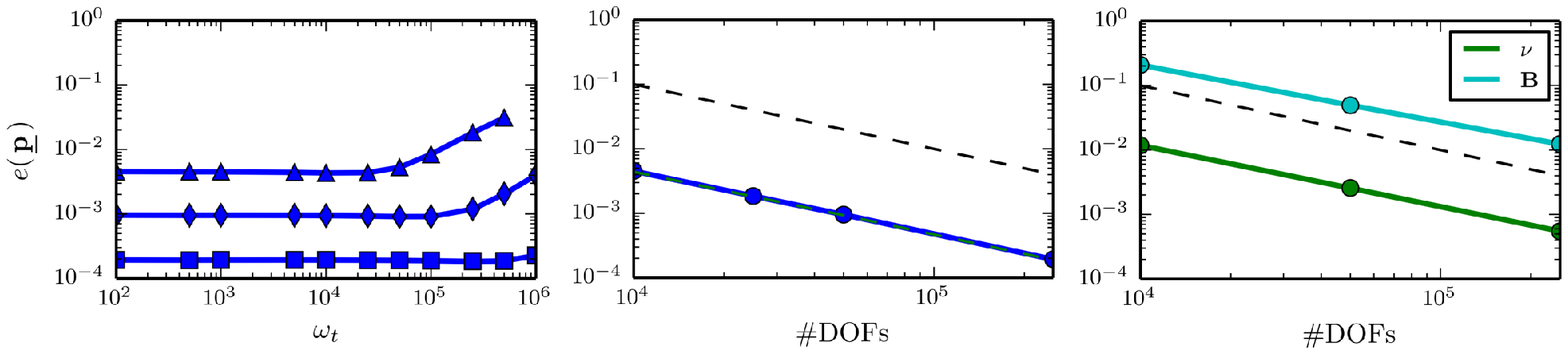}
\end{center}
\caption{Numeric convergence of $\pb$ on ellipsoid: (left) $L^2$ error $e(\ext{\pb})$ as function of $\omega_t$. Markers indicate surface discretization with 10k (triangles), 50k (diamonds) and 250k (squares) DOFs. (middle) $L^2$ error measure with analytical geometric quantities (dashed green line) and numerically approximated surface quantities (blue solid line) w.r.t. number of DOFs for fixed $\omega_t = 1000$. Black dashed line indicates linear rate of convergence. (right) Convergence of $L^2$ error for numerical approximation of surface quantities w.r.t. number of DOFs for fixed $\omega_t = 1000$. (colors online)}
\label{fig:numericConvergencep}
\end{figure}

The inner product and the covariant derivative expressed in Cartesian coordinates read
\begin{align}
\langle \ext{\pb}, \ext{\psib}  \rangle_1 = \int_{\manifold} \tensor{\ProjSurf}{^{I}_{J}} \ext{p}_I \ext{\psi}^J \;\mathrm{d}\manifold, \qquad
\nabla \pb = \ProjSurf[\extNabla\ext{\pb}] + \left( \normal \iP{1} \ext{\pb} \right) \shapeOperator, \qquad
[\nabla \pb]_{IJ} = \tensor{\ProjSurf}{_{I}^{K}}\tensor{\ProjSurf}{_{J}^{L}}\extNabla_L \ext{p}_K + \ext{p}_L\nu^L B_{IJ}.
\end{align}
As penalty term we use $p(\ext{\pb}) = \omega_t \left(\normal \outP \normal \right) \inP \ext{\pb}$, with $\omega_t > 0$. The weak formulation thus reads
\begin{align}
&\int_{\manifold} \left( \ProjSurf \iP{} [\extNabla \ext{\pb}] \right) : \left( [\extNabla \ext{\psib}] \iP{2} \ProjSurf  \right) \;\mathrm{d}\manifold \nonumber \\
& +\int_{\manifold} \left( \shapeOperator : [\extNabla \ext{\pb}]\right)\left( \normal \iP{} \ext{\psib} \right)  + \left( \normal \iP{} \ext{\pb} \right)\left( \shapeOperator : [\extNabla \ext{\psib}]\right)  \;\mathrm{d}\manifold \nonumber \\
&+ \int_{\manifold}\left( \meanCurvature^2 - 2 \gaussianCurvature \right) \left( \normal \iP{} \ext{\pb} \right) \left( \normal \iP{} \ext{\psib} \right) \;\mathrm{d}\manifold \nonumber \\
+ & \int_{\manifold} \ext{\pb} \iP{}\ProjSurf \iP{}  \ext{\psib} \;\mathrm{d}\manifold +  \int_{\manifold} \omega_t \left(\ext{\pb} \iP{} \normal\right) \left( \ext{\psib} \iP{} \normal \right) \; \mathrm{d} \manifold = \int_{\manifold} \fb \iP{} \ProjSurf \iP{} \ext{\psib} \;\mathrm{d}\manifold \quad \forall \ext{\psib} \in [H^1(\manifold)]^{3,1}
\end{align}
where we have used several geometric identities given in \ref{sec:geometricIdentities}. We now consider a surface triangulation $\manifold_h$ and construct a surface finite element discretization for each component. We consider the finite element space $V(\manifold_h)$ of piecewise linear Lagrange elements. The problem than reads: Find $\ext{\pb} = (\ext{p}_1,\ext{p}_2,\ext{p}_3)  \in [V(\manifold_h)]^{3,1}$ s.t.
\begin{align}
\label{eq:pVectorFEM}
 &\quad\int_{\manifold_h} \tensor{\ProjSurf}{_{I}^{K}}\tensor{\ProjSurf}{_{J}^{L}}\extNabla_L \ext{p}_K \extNabla^J \ext{\psi}^I   + \nu_L B^{IJ} \extNabla_J \ext{p}_I \ext{\psi}^L + \nu^L B_{IJ}  \ext{p}_L \extNabla^J \ext{\psi}^I + \left( \meanCurvature^2 -2\gaussianCurvature \right) \nu^I\nu_J \ext{p}_I \ext{\psi}^J  \;\mathrm{d}\manifold_h \nonumber \\
 &  + \int_{\manifold_h} \tensor{\ProjSurf}{^{I}_{J}} \ext{p}_I \ext{\psi}^J \;\mathrm{d}\manifold_h + \int_{\manifold_h} \omega_t \nu^I\nu_J \ext{p}_I  \ext{\psi}^J  \; \mathrm{d} \manifold_h 
= \int_{\manifold_h} \tensor{\ProjSurf}{^{I}_{J}} f_I \ext{\psi}^J \;\mathrm{d}\manifold_h \quad \forall \ext{\psib}=(\ext{\psi}^1, \ext{\psi}^2, \ext{\psi}^3) \in [V(\manifold_h)]^{3,1}
\end{align}
where $\extNabla$, $\shapeOperator$, $\meanCurvature$, $\gaussianCurvature$, $\normal$ and $\ProjSurf$ have to be understood with respect to $\manifold_h$. We consider both the analytically available geometric quantities and numerically computed once from the surface triangulation, see \cite{Nitschkeetal_JFM_2012} for the used numerical approach to compute $\normal$ and $\shapeOperator$. The resulting linear system is assembled and solved in the finite element toolbox AMDiS \cite{Veyetal_CVS_2007,Witkowskietal_ACM_2015} using a BiCGstab({\it{l}}) solver. As error measure we use the componentwise averaged $L^2$ norm
\begin{align}
e(\extTF{\pb})=\frac{1}{2}\left(\int_{\manifold_h} \| \extTF{\pb} - \pb^* \|^2 \mathrm{d}\manifold_h \right)^{1/2}
\end{align}
where $\| \pb \|$ indicates the Forbenius norm in $\R^{3,1}$. Convergence properties regarding the number of degrees of freedom (DOFs) and the strength of the penalty parameter $\omega_t$ are shown in figure \ref{fig:numericConvergencep}.
The results indicate almost no dependency on $\omega_t$, which reflect the behavior of the analytical solution, up to values where the penalty term begins to dominate the numerical solution behavior. We further see linear convergence in the geometric properties and also a linear convergence behavior in $e(\pb)$, which does not depend on the used approximation of the geometric terms. As the geometry is approximated by a surface triangulation better than linear convergence properties could not be expected. Please note, assuming a mesh with uniform element size $h$, linear convergence in $\#DOFs$ corresponds to quadratic convergence in $h$, which is shown to be optimal for a surface finite element approach with piecewise linear Lagrange elements for a scalar-valued surface Laplace problem \cite{Dziuk,Demlow_SIAMJNA_2009}.

\subsection{Tensor-valued Helmholtz equation} 
\begin{figure}[ht]
\begin{center}
\includegraphics[height=0.3\textwidth]{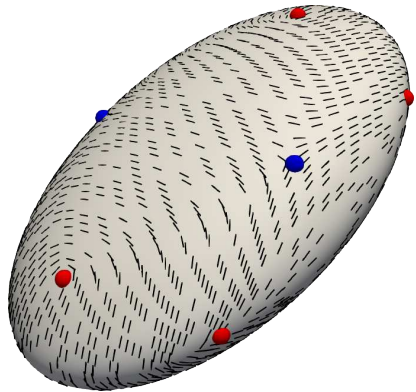}
\end{center}
\caption{\textbf{Principal eigenvectors of analytical solution $\qb^*$ with 8 nodes (topological charge $+ \frac{1}{2}$, red) and 4 wedges (topological charge $- \frac{1}{2}$, blue). (color online).}}
\label{fig:tensorHelmholtz}
\end{figure}

We consider Q-tensors $\qb$ which are symmetric and traceless 2-tensors and formulate a surface  tensor-valued Helmholtz problem. There are different options how to enforce the traceless property in the solution procedure. We demonstrate an approach where the property is accounted for in
solution subspaces  of $[H^1(\manifold)]^{3,2}$ and $[V(\manifold_h)]^{3,2}$. However, we are faced with different notions of traces in the tangential and the embedding space.
For 2-tensors $ \ext{\qb} $ and surface  2-tensors $ \tilde{\qb}=\ProjSurf[\ext{\qb}] $, not necessarily traceless, holds 
\begin{align}
\trace(\tilde{\qb})=(\tilde{\qb},\gb)_{2,\manifold} = \Gb : \ext{\qb} - \normal \inP \ext{\qb} \inP \normal = \traceR(\ext{\qb}) - \normal \inP \ext{\qb} \inP \normal
\end{align}
We can account for this issue by considering the extension $\qb = \tilde{\qb} + \frac{1}{2} [\ext{\qb}]_{\xi\xi}\gb$. 
Therefore $ \trace(\qb)= 0 $ iff $ \traceR(\ext{\qb})=0 $ is valid.
Note that $ \qb $ is not simply the tangential projected $ \ext{\qb} $ rather then the result of an orthogonal projection into surface Q-tensors.
This requires to redo the calculations leading to \eqref{eqn:exactCovarinatDerivative}. 
With
\begin{align}
[\extNabla\ext{\qb}]_{\xi\xi k}  & = \partial_k \ext{q}_{\xi\xi} - \ChristE{L}{k\xi}\ext{q}_{L\xi} - \ChristE{L}{k\xi}\ext{q}_{\xi L} = [\nabla \ext{q}_{\xi\xi}]_k + 2 \tensor{B}{_{k}^{l}}\ext{q}_{l\xi} \formComma
\end{align}
thin film coordinates and metric compatibility leads to
\begin{align}
	[\nabla\qb]_{ijk} &= [\nabla\tilde{\qb}]_{ijk} + \frac{1}{2}g_{ij}[\nabla \ext{q}_{\xi\xi}]_k
		= [\extNabla\ext{\qb}]_{ijk} + \frac{1}{2}g_{ij}[\extNabla\ext{\qb}]_{\xi\xi k}
		  + \tensor{B}{_{ki}}\tensor{\ext{q}}{_{\xi j}} + \tensor{B}{_{kj}}\tensor{\ext{q}}{_{i\xi}}
		  - g_{ij}\tensor{B}{_{k}^{l}}\ext{q}_{l\xi} \formPeriod
\end{align}
This allows to obtain the desired Cartesian description of the covariant derivative
\begin{align}
\!\!\! [\nabla\qb]_{IJK} = & \tensor{\ProjSurf}{_{I}^{M}}\tensor{\ProjSurf}{_{J}^{N}}\tensor{\ProjSurf}{_{K}^{O}}\extNabla_O\ext{q}_{MN} \nonumber \\
& + \tensor{B}{_{KI}} \nu_{M} \tensor{\ext{q}}{^{M}_{N}} \tensor{\ProjSurf}{^{N}_{J}} + \tensor{B}{_{KJ}} \tensor{\ProjSurf}{_{I}^{M}} \tensor{\ext{q}}{_{M}^{N}}  \nu_{N} \nonumber \\
& + \frac{1}{2}\tensor{\ProjSurf}{_{IJ}}\tensor{\ProjSurf}{_{K}^{N}} \nu^L\nu^M \extNabla_N\ext{q}_{LM} -\tensor{\ProjSurf}{_{IJ}}\tensor{B}{_{K}^{L}}\nu^M \ext{q}_{LM} \nonumber \\
[\nabla\qb] = & \ProjSurf[\extNabla\ext{\qb}] +  \shapeOperator \oP{1,3}\left( \normal \iP{} \ext{\qb} \iP{2} \ProjSurf \right)  + \shapeOperator \oP{2,3} \left( \ProjSurf \iP{} \ext{\qb} \iP{2} \normal \right) 
+ \frac{1}{2} \ProjSurf \oP{} \ProjSurf\left[ \normal\iP{}(\extNabla\ext{\qb})\iP{2}\normal \right]  
- \ProjSurf \oP{} \left( \shapeOperator \iP{} \ext{\qb} \iP{2} \normal \right) \formPeriod \label{eq:covdertensor}
\end{align}
Together with the inner product  
\begin{align}
\langle \ext{\qb},\ext{\psib} \rangle_{2,\manifold} = \int_{\manifold}  \tensor{\ProjSurf}{^{I}_{K}}  \tensor{\ProjSurf}{^{J}_{L}}\ext{q}_{IJ} \ext{\psi}^{LK} \; \mathrm{d}\manifold
\end{align}
we have all tools available to formulate the weak form
\begin{align}
\tilde{a}(\ext{\qb},\ext{\psib}) = l(\ext{\psib}) \quad \forall \ext{\qb},\ext{\psib} 
   \in Q^{1,3}(\manifold) = \left\{ \ext{\qb} \in [H^1(\manifold)]^{3,2}\ \middle|\ \ext{\qb} = \ext{\qb}^T \text{ and } \traceR(\ext{\qb})=0 \right\}
\end{align}
with $\tilde{a}(\ext{\qb},\ext{\psib}) = a(\ext{\qb},\ext{\psib}) + \langle p(\ext{\tb}), \ext{\psib}  \rangle$, a modified definition for $a(\cdot, \cdot)$ according to \eqref{eq:covdertensor} and a penalty term \wrt\ complementary tangential Q-tensor projection as 
$p(\ext{\qb}) = \omega_t (\normal \outP  (\normal \inP \ext{\qb}) - \frac{1}{4}(\normal \iP{} \ext{\qb} \iP{2}  \normal)\normal \outP  \normal)$, $\manifold$ again an ellipsoid with major axis $A=1$, $B=0.5$ and $C=1.5$, $\fb = -\Div \left(  \nabla \qb^* \right) + \qb^*$, with $\qb^* = \pb^{*} \oP{} \pb^{*} - \frac{1}{2} \gb$ and $\pb^* = \Rot(xyz)$, see figure  \ref{fig:tensorHelmholtz}. The solution contains various defects with topological charge (winding number) $\pm \frac{1}{2}$, while the sum of these charges (winding numbers) is again two. 

\begin{figure}[ht]
\begin{center}
\includegraphics[height=0.2\textwidth]{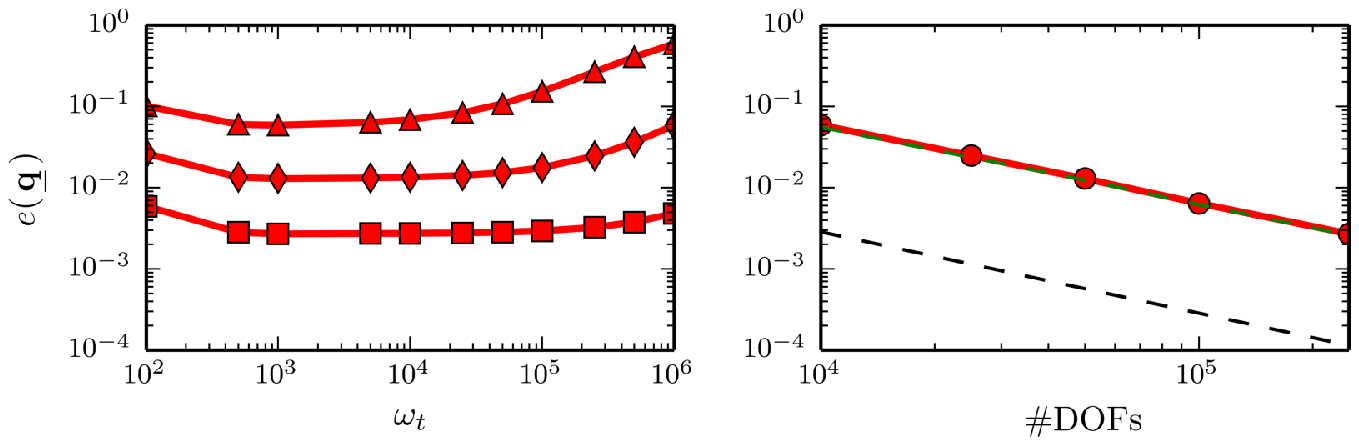}
\end{center}
\caption{\textbf{Numeric convergence test of $\qb$ on ellipsoid:} (left) $L^2$ error measure $e(\ext{\qb})$ as function of $\omega_t$. Markers indicate surface discretization with 10k (triangles), 50k (diamonds) and 250k (squares) DOFs. (right) $L^2$ error measure $e(\ext{\qb})$ with analytical geometric quantities (green dashed line) and numerically approximated surface quantities (red solid line) w.r.t. number of DOFs for fixed $\omega_t = 1000$. Black dashed line indicates linear rate of convergence. (colors online)}
\label{fig:numericConvergenceq}
\end{figure}

The weak formulation reads
\begin{align}
\int_{\manifold} [\nabla\qb]_{IJK}[\nabla\psib]^{IJK} \;\mathrm{d}\manifold + \int_{\manifold} \tensor{\ProjSurf}{^{I}_{K}}  \tensor{\ProjSurf}{^{J}_{L}}\ext{q}_{IJ} \ext{\psi}^{LK} \;\mathrm{d}\manifold + \int_{\manifold} [p(\ext{\pb})]_{IK} \ext{\psi}^{IK} \; \mathrm{d} \manifold 
 =  \int_{\manifold} \tensor{\ProjSurf}{^{I}_{K}}  \tensor{\ProjSurf}{^{J}_{L}}f_{IJ} \ext{\psi}^{LK} \;\mathrm{d}\manifold 
\end{align}
$\forall \ext{\psib} \in Q^{1,3}(\manifold)$ with
\begin{align}\label{bilin}
	& \int_{\manifold} [\nabla\qb]_{IJK}[\nabla\psib]^{IJK} \;\mathrm{d}\manifold \nonumber \\
		&= \int_{\surf} \left[\tensor{\ProjSurf}{^{M}_{S}} \tensor{\ProjSurf}{^{N}_{T}} + 1/2\left(  \nu^M \nu^N \nu_{S} \nu_{T} + \ProjSurf^{MN} \nu_{S} \nu_{T} +\nu^{M} \nu^{N} \ProjSurf_{ST} \right) \right]\tensor{\ProjSurf}{^{K}_{L}}
		\extNabla_K \ext{q}_{MN} \extNabla^{L} \ext{\psi}^{ST}  \nonumber \\
		&\quad + \left[2\tensor{B}{^{KN}} \tensor{\ProjSurf}{^{M}_{S}} \nu_{T} - \tensor{B}{^{K}_{T}} \ProjSurf^{MN}\nu_{S}\right]
		       \extNabla_K  \ext{q}_{MN} \ext{\psi}^{ST} \\
		&\quad + \left[2\tensor{B}{^{K}_{T}} \tensor{\ProjSurf}{^{M}_{S}} \nu^{N} -B^{KN}  \ProjSurf_{ST} \nu^{M}\right]
		       \ext{q}_{MN} \extNabla_K \ext{\psi}^{ST} 
		    +2\left(\meanCurvature^2 - 2\gaussianCurvature \right)\tensor{\ProjSurf}{^{N}_{T}}\nu^{M}  \nu_{S}
		     \ext{q}_{MN} \ext{\psi}^{ST} \; \mathrm{d}\manifold \notag
\end{align}
where the identities in \ref{sec:geometricIdentities} have been used to simplify the expressions. The corresponding finite element formulation reads: Find $\ext{\qb} = (\ext{q}_1, \ext{q}_2, \ext{q}_3, \ext{q}_4, \ext{q}_5) \in [V(\manifold_h)]^{5,1}$ s.t.
\begin{align}
\!\!\!\!& \int_{\manifold_h} [\nabla\qb]_{IJK}[\nabla\psib]^{IJK} \;\mathrm{d}\manifold_h  + \int_{\manifold_h} \tensor{\ProjSurf}{^{I}_{K}}  \tensor{\ProjSurf}{^{J}_{L}}\ext{q}_{IJ} \ext{\psi}^{LK} \;\mathrm{d}\manifold_h 
+ \int_{\manifold_h} [p(\ext{\pb})]_{IK} \ext{\psi}^{IK} \; \mathrm{d} \manifold_h \nonumber \\
&= \int_{\manifold_h} \tensor{\ProjSurf}{^{I}_{K}}  \tensor{\ProjSurf}{^{J}_{L}}f_{IJ} \ext{\psi}^{LK} \;\mathrm{d}\manifold_h 
\end{align}
$\forall \ext{\psib} = (\ext{\psi}^1, \ext{\psi}^2, \ext{\psi}^3, \ext{\psi}^4, \ext{\psi}^5) \in [V(\manifold_h)]^{5,1}$ used in a 
3x3 matrix proxy notation for Cartesian coordinate system, \eg 
\begin{align}
\label{eq:componentAnsatzQ}
\ext{\qb} = \begin{bmatrix}
\ext{q}_1 & \ext{q}_2 & \ext{q}_3  \\
\ext{q}_2 & \ext{q}_4 & \ext{q}_5  \\
\ext{q}_3 & \ext{q}_5 &  -\ext{q}_1-\ext{q}_4
\end{bmatrix} \quad \mbox{and}\quad
\ext{\psib} = \begin{bmatrix}
\ext{\psi}^1 & \ext{\psi}^2 & \ext{\psi}^3  \\
\ext{\psi}^2 & \ext{\psi}^4 & \ext{\psi}^5  \\
\ext{\psi}^3 & \ext{\psi}^5 &  -\ext{\psi}^1-\ext{\psi}^4
\end{bmatrix}\formComma
\end{align}
which only considers the symmetric and traceless subspace of $[V(\manifold_h)]^{3,2}$.
We use the same tools as in the vector-valued case and show convergence properties regarding the number of degrees of freedom (DOFs) and the strength of the penalty parameter $\omega_t$ in figure \ref{fig:numericConvergenceq}, with average componentwise error measure 
\begin{align}
e(\extTF{\qb})=\frac{1}{4}\left(\int_{\manifold_h} \| \extTF{\qb} - \qb^* \|^2 \mathrm{d}\manifold_h \right)^{1/2}
\end{align}
where $\| \qb \|$ indicates the Forbenius norm in $\R^{3,2}$.
The results indicate almost the same behavior as for the vector case. There is almost no dependency on $\omega_t$, up to values where the penalty term begins to dominate the numerical solution behavior. Only for $\omega_t < 1000$ a small dependency can be seen. We further see the optimal linear convergence behavior in $e(\ext{\qb})$, which again does not depend on the used approximation of the geometric terms. However, the overall error is larger than in the vector case, which is also expected due to the more complex problem. 

\subsection{General vector- and tensor-valued surface PDEs}
As seen in the previous paragraph the complexity of the problem not only increases due to the higher tensorial degree, also the incorporation of tensor properties can led to several additional terms. Even if optimal convergence properties are only experimentally shown for two ``simple" problems, we expect this behavior to hold also for other surface PDEs. All tools are provided, see also \ref{sec:AppendixIntByParts} - \ref{sec:AppendixDifferentialOperators}, to rewrite any vector- or tensor-valued surface PDE in Cartesian coordinates and to solve them by applying a surface finite element method, or any other appropriate method to solve scaler-valued surface PDEs, to each component. Also the possibility to account for specific tensor properties is, at least exemplarily, shown. This opens a huge field of applications and the possibility to numerically analyze the solutions. However, both is beyond the scope of this paper. 

\subsection{Error sensitivities}
\begin{figure}[ht]
\begin{center}
\includegraphics[height=0.2\textwidth]{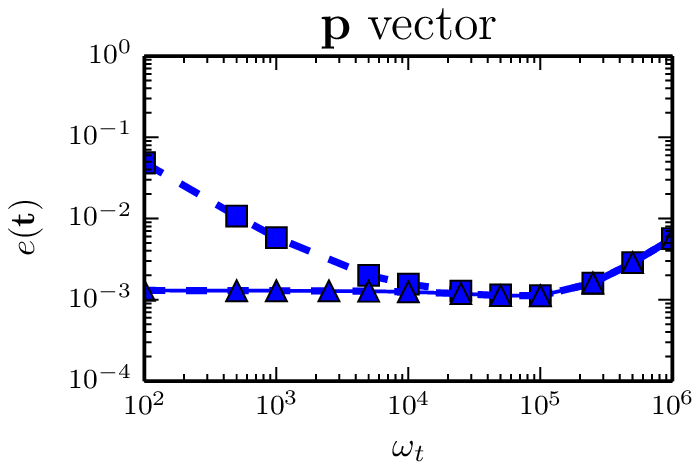}
\includegraphics[height=0.2\textwidth]{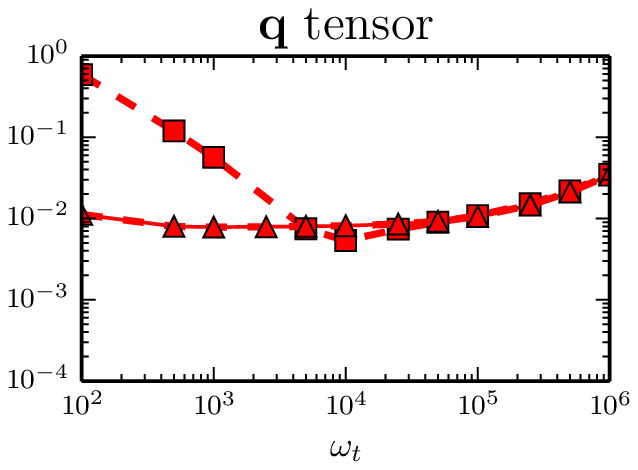}
\includegraphics[height=0.2\textwidth]{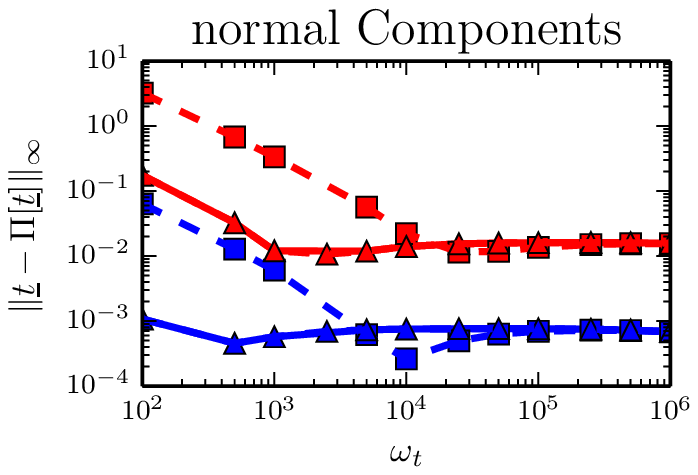}
\end{center}
\caption{\textbf{Numeric convergence for approximated operators in considered Helmholtz problems:} Normalized $L^1$ error measure for $\pb$ (left) and $\qb$ (middle) for exact (solid line), approximated scalar product (triangles) and approximated covariant derivative (squares). (right) The maximum norm of the normal components for $\pb$ (blue) and $\qb$ (red) is shown for the considered approximations (symbols as before). $\manifold_h$ is discretized with 50k vertices.}
\label{fig:numericConvergenceApproximated}
\end{figure}

The numerical approach is not only complex, the assembly process is also computationally demanding. We therefore consider two possibilities to reduce the complexity and the assembly workload significantly, which can be achieved by considering $\nabla \tb \approx \ProjSurf[\extNabla\ext{\tb}] $ and $\langle \tb, \qb \rangle \approx \int_{\manifold} \ext{\tb}:\ext{\qb} \,\mathrm{d}\manifold$
for $ \ext{\tb}\approx\ProjSurf[\ext{\tb}] $ and $ \ext{\qb}\approx\ProjSurf[\ext{\qb}] $. The results of these approximations on the considered Helmholtz problems are shown in figure \ref{fig:numericConvergenceApproximated} where the normalized $L^1$ error is used.
\begin{align}
e(\extTF{\tb})=\frac{\int_{\manifold_h} \| \extTF{\tb} - \tb^* \| \mathrm{d}\manifold_h}{(n-1)^d \int_{\manifold_h} \| \tb^* \| \mathrm{d}\manifold_h }
\end{align}
Approximating the scalar product, while keeping the exact covariant derivative, does not significantly increase the error, regardless of the choice of $\omega_t$. Contrary the use of approximated derivatives exhibits a strong dependance on $\omega_t$ such that only for high values the same error level can be reached and determining an appropriate value for $\omega_t$ might be hard in more complex applications. These results underline the importance of the additional coupling term in \eqref{eqn:exactCovarinatDerivative} and demonstrate the advantage to previous methods \cite{Nestleretal_JNS_2018,Reutheretal_PF_2018,Hansboetal_arXiv_2017}. 

\section{Application example}
\label{sec:4}
In \cite{Nitschkeetal_PRSA_2018} a surface Landau-de Gennes model for nematic liquid crystals, whose molecular orientation is subjected to a tangential anchoring on a curved surface, was derived and numerically solved using the described approach. The model in its one-constant approximation reads
\begin{align}\label{eq:evoliEquation}
  \partial_{t}\qb - \Delta^{dG}\qb  + \frac{1}{2} \left( \meanCurvature^2 - 2 \gaussianCurvature \right) \qb + \omega (1 -  2 \trace\qb^{2}) \qb &= 0
\end{align}
on $\manifold \times [0,t_{\textup{end}}]$ with the div-Grad (Bochner) Laplacian $\Delta^{dG} \qb = \Div (\nabla \qb)$. We consider the following parameters $L_1 = L_2 = - L_3 = 1$, $L_6 = 0$, $a = c = 1$, $b = 0$, $\omega = 100$ and $\beta = 0$, see \cite{Nitschkeetal_PRSA_2018} for the full model. The equation results as an $L^2$-gradient flow of the surface Landau-de Gennes energy ${\cal{F}} = {\cal{F}}_{el} + \omega {\cal{F}}_{bulk}$ with an elastic and an entropic part, with $L_i$ the elastic parameters, and $a$, $b$, $c$ the entropic parameters. The quantity of interest is the surface $Q$-tensor $\qb$. The model thus can be seen as an extension of the considered tensor-valued Helmholtz problem and we adapt the used discretization within each timestep. As $\manifold_h$ we consider a sufficiently complex manifold, a slightly smoothed and refined surface triangulation of the Stanford bunny, and as initial condition $\qb(\xb,0) = \qb_0(\xb)$ random noise. 

\begin{figure}[ht]
\begin{center}
\includegraphics[width=0.8\textwidth]{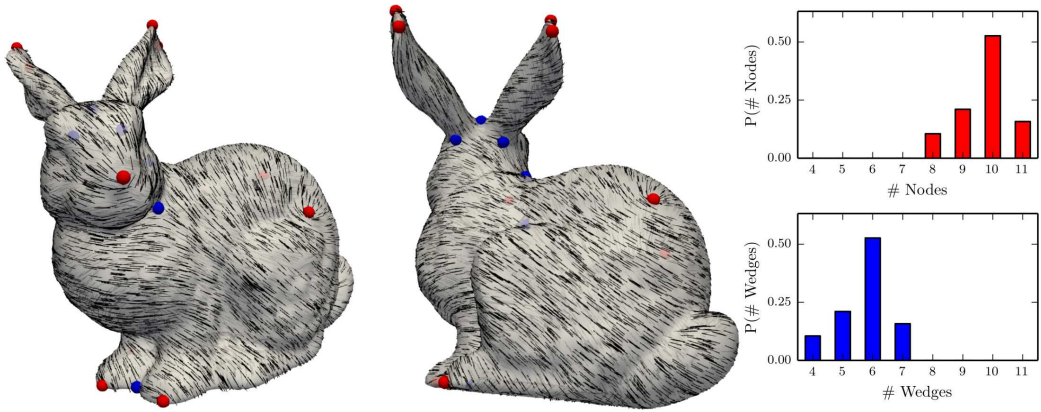}
\end{center}
\caption{\textbf{Nematic order and defect configuration on the Stanford bunny:} Relaxation from noise to energy minima was performed 20 times. Typical realization with 10 nodes (red) and 6 wedge (blue) defects in front (left) and back (middle) view. (right) Statistical distribution of defect configuration realizations w.r.t. number of node and wedge defects. (colors online) }
\label{fig:stanfordBunny}
\end{figure}

For the discretization in time let the time interval $[0,t_{\textup{end}}]$ be divided into a sequence of discrete times $0 < t_0 < t_1 < ...$ with time step width $\tau^m = t^{m} - t^{m-1}$. Thereby, the superscript denotes the timestep number and $\qb^m(\xb) = \qb(\xb, t^m)$. The time derivative is approximated by a standard difference quotient, we define the discrete time derivatives $d_{\tau^m}^{\qb} := \frac{1}{\tau^m}\left( \qb^{m} - \qb^{m-1} \right)$. Thus, we get a time-discrete version of \eqref{eq:evoliEquation}
\begin{align}\label{eq:evoliEquationdis}
 d_{\tau^m}^{\qb}  - \Delta^{dG}\qb^m  + \frac{1}{2} \left( \meanCurvature^2 - 2 \gaussianCurvature \right) \qb^m + \omega \left(1 -  2 \trace(\qb^m)^{2}\right) \qb^m &= 0
\end{align}
with a non-linear term, which will be treated by a Newton iteration in each time step. Adapting the weak formulation from the tensor-valued surface Helmholtz problem we get for $\ext{\qb}^m$
\begin{align}
&\int_{\manifold} \tensor{\ProjSurf}{^{I}_{K}}  \tensor{\ProjSurf}{^{J}_{L}}\ext{d_{\tau^m}^{\qb}}_{IJ} \ext{\psi}^{LK} \;\mathrm{d}\manifold + \int_{\manifold} [\nabla\qb^m]_{IJK}[\nabla\psib]^{IJK} \;\mathrm{d}\manifold + \int_{\manifold} \tensor{\ProjSurf}{^{I}_{K}}  \tensor{\ProjSurf}{^{J}_{L}} \frac{1}{2} \left( \meanCurvature^2 - 2 \gaussianCurvature \right) \ext{q}^m_{IJ} \ext{\psi}^{LK} \;\mathrm{d}\manifold \nonumber \\
&\quad+ \int_{\manifold} p(\ext{\qb}^m)_{IK}  \ext{\psi}^{IK} \; \mathrm{d} \manifold 
+  \int_{\manifold} \tensor{\ProjSurf}{^{I}_{K}}  \tensor{\ProjSurf}{^{J}_{L}} \omega\left(1 -2 \trace(\ext{\qb}^{m})\right)^2 \ext{q}^m_{IJ} \ext{\psi}^{LK} \;\mathrm{d}\manifold = 0 \qquad \forall \ext{\psib} \in [H^1(\manifold)]^{3,2}
\end{align}
with $ \int_{\manifold} [\nabla\qb^m]_{IJK}[\nabla\psib]^{IJK} \;\mathrm{d}\manifold$ defined in \eqref{bilin}. For better readability we omit the formulation of the Newton iteration and keep the non-linear term. The corresponding finite element formulation is obtained by replacing $[H^1(\manifold)]^{3,2}$ by $[V(\manifold)]^{3,2}$. We consider the same component function ansatz for $\ext{\qb}$ and $\ext{\psib}$ as defined in \eqref{eq:componentAnsatzQ} and again piecewise linear Lagrange elements. The resulting linear system is solved in parallel using a preconditioned BiCGstab({\it{l}}) solver. We use $\leq 10$ Newton iterations in each time step. Results for $\omega_t = 1000$ are shown in figure \ref{fig:stanfordBunny}. 
The configuration of four node defects on the ears combined with three wedge defects between the ears is invariant in the group of realizations with 10 nodes and 6 wedges. The node and wedge defect distribution along nose, back and feet varies. All realizations fulfill the Poincare-Hopf theorem, the sum of the topological charges (winding numbers) is always two. The lowest energy ${\cal{F}}$ among the 20 realizations was found for a configuration with 9 nodes and 5 wedges. All realizations have in common that nodes are located at points of positive Gaussian curvature (ears, nose, feet, back), while wedges are located at points of negative Gaussian curvature (between ears, between feet). However, a detailed analysis of these connection goes beyond the scope of this paper.

\section{Conclusion}
\label{sec:5}
We proposed a general finite element method to solve tangential vector- and tensor-valued surface PDEs. The key idea is the reformulation of the problem in Cartesian coordinates. This can be done by considering an extended space including normal components. We derive the formulations for the first order operators $\nabla$, $\Div$ and $\Rot$ and obtain various terms which couple the normal component with geometric properties of the manifold. The complexity of these terms increases with the tensorial degree. The proposed formulations are exact and invariant under changes in normal direction. In order to make the solution unique a penalization is considered which forces the normal components to be zero. The considered test cases show that the solution is almost independent of the penalization parameter and converges with optimal order. The reformulation allows to solve the problems in a componentwise fashion using established tools for scalar-valued surface PDEs. We here only considered surface finite elements, but any other approach to solve scalar-valued surface PDEs would be appropriate as well. The huge flexibility for these methods in terms of topology, shape and description of it, thus transfers to tangential vector- and tensor-valued surface PDEs. The approach can also be naturally extended to evolving surfaces, see \cite{Nestleretal_arXiv_2018} for a first realization again using surface finite elements for each component. For a general discussion on transport of vector- and tensor-quantities on evolving surfaces we refer to \cite{Nitschkeetal_arXiv_2018}. \\

\noindent
{\bf Acknowledgments}: This work was financially supported by the German Research Foundation (DFG) through project Vo899-19. We used computing resources provided by J\"ulich Supercomputing Centre within project HDR06.

\appendix

\section{Integration by parts for tensor fields}
\label{sec:AppendixIntByParts}
Using $\conormal$ as conormal of the boundary $\partial \manifold$, the
divergence theorem
\begin{align}
    \int_{\manifold} \Div\pb\, \mathrm{d}\manifold 
        &= \int_{\partial\manifold} \left( \pb, \conormal \right)_{1,\manifold} \mathrm{d}\partial\manifold
\end{align} 
holds also in Riemannian manifolds for vector fields $ \pb\in\tensorField{1}{\manifold} $, 
see \cite{Unal1995,Andrews2011}.
This is a consequence of Cartan's magic formula, Stokes theorem (see \cite{abraham2012manifolds}) and 
the stipulation that $ \mathrm{d}\partial\manifold = \mathrm{d}\manifold(\conormal, \hdots) $, which determine the orientation
of $ \partial\manifold $. 

Concerning the partial integration for tensorial fields on manifolds, let $\tb \in \tensorField{d+1}{\manifold}$ and $\psib \in \tensorField{d}{\manifold}$, 
we yield 
\begin{align}
\int_{\manifold} \left(\Div \tb, \psib  \right)_{d,\manifold} \, \mathrm{d}\manifold 
    &= \int_{\manifold} \nabla_k t^{i_1\hdots i_d k} \psi_{i_1\hdots i_d}\, \mathrm{d}\manifold
     = \int_{\manifold} \nabla_k \left( t^{i_1\hdots i_d k} \psi_{i_1\hdots i_d}  \right) 
                       - t^{i_1\hdots i_d k}\nabla_k\psi_{i_1\hdots i_d}\, \mathrm{d}\manifold \notag\\
    &= -\int_{\manifold} \left(\tb, \nabla\psib  \right)_{d+1,\manifold} \, \mathrm{d}\manifold
       + \int_{\partial\manifold}  \left( \tb\iP{d+1}\conormal, \psib \right)_{d,\manifold} \mathrm{d}\partial\manifold
\end{align}
as a consequence of Leibniz rule and the divergence theorem above for $ p^k = t^{i_1\hdots i_d k} \psi_{i_1\hdots i_d} $.

\section{Componentwise description of boundary terms}
\label{sec:AppendixComponentwiseBoundaryConditions}
\subsection{Cartesian description of boundary terms for 1-tensor (vector) fields}
Applying integration by parts to obtain the weak formulation of the div-Grad Laplacian yields for vector fields
\begin{align}
\int_{\manifold} \left(-\Div \nabla \pb, \psib  \right)_{1,\manifold} \, \mathrm{d}\manifold 
    & = -\int_{\partial\manifold} \left( \nabla_{\conormal}\pb, \psib \right)_{1,\manifold} \; \mathrm{d}\partial\manifold + \int_{\manifold} \left(\nabla \pb, \nabla \psib  \right)_{2,\manifold} \, \mathrm{d}\manifold \formComma 
\end{align}
where $ \nabla_{\conormal}\pb = \nabla\pb\iP{2}\conormal $ is the conormal derivative of $ \pb $. 
In Cartesian description we yield two contributions for the boundary term
\begin{align}
     &\int_{\partial\manifold} \left( \nabla_{\conormal}\pb, \psib \right)_{1,\manifold} \; \mathrm{d}\partial\manifold \nonumber \\
        =& \int_{\partial\manifold} \ProjSurf^{IJ} \ext{\psi}_J n^{K} \left( \extNabla_{K} \ext{p}_{I} 
                                    + B_{KI}\ext{p}_L \nu^{L}  \right)\; \mathrm{d}\partial\manifold \nonumber \\
         =& \int_{\partial\manifold} \left( \extNabla_{\conormal}\ext{\pb}, \ext{\psib} \right)_{1,\manifold} 
                                    + \left( \ext{\pb}\iP{}\normal \right)\left( \shapeOperator\iP{}\conormal, \ext{\psib} \right)\; \mathrm{d}\partial\manifold
\end{align}
where the first one vanish if homogeneuos Neumann condition holds and the second one if and only if $\ext{\pb}|_{\partial \manifold}$ is strictly tangential
or the curvature vanish at $ \partial\manifold $ in direction of the conormal $ \conormal $, generally for all $ \ext{\psib} $ 
and $ \partial\manifold\neq\emptyset $.

\subsection{Cartesian description of boundary terms for 2-tensor fields}
Integration by parts for the weak formulation of div-Grad Laplacian yields
\begin{align}
\int_{\manifold} \left(-\Div \nabla \qb, \psib  \right)_{2,\manifold} \, \mathrm{d}\manifold 
    & = -\int_{\partial\manifold} \left( \nabla_{\conormal}\qb, \psib \right)_{2,\manifold} \; \mathrm{d}\partial\manifold 
        + \int_{\manifold} \left(\nabla \qb, \nabla \psib  \right)_{3,\manifold} \, \mathrm{d}\manifold \formComma 
\end{align}
where $ \nabla_{\conormal}\qb = \nabla\qb\iP{3}\conormal $ is the conormal derivative of $ \qb $.
In Cartesian coordinates the boundary term is expressed by
\begin{align}
    \int_{\partial\manifold} \left( \nabla_{\conormal}\qb, \psib \right)_{2,\manifold} \; \mathrm{d}\partial\manifold
        &= \int_{\partial\manifold} \ProjSurf^{I_1 J_1} \ProjSurf^{I_2 J_2} \ext{\psi}_{J_1 J_2} n^{K}
                \left( \extNabla_{K}\ext{q}_{I_1 I_2} + \nu^{L}\left( B_{K I_1} \ext{q}_{L I_2} +  B_{K I_2} \ext{q}_{I_1 L}\right) \right)
            \; \mathrm{d}\partial\manifold \notag \\
        &= \int_{\partial\manifold} \left( \extNabla_{\conormal}\ext{\qb}, \ext{\psib} \right)_{2,\manifold}
           + \left( \left( \shapeOperator\iP{}\conormal \right) \oP{} \left( \normal \iP{} \ext{\qb} \iP{2} \boldsymbol{\ProjSurf} \right)
                   + \left( \boldsymbol{\ProjSurf} \iP{} \ext{\qb} \iP{2} \normal\right) \oP{} \left( \shapeOperator\iP{}\conormal \right)           
                , \ext{\psib} \right) \; \mathrm{d}\partial\manifold
\end{align}
Therefore a homogeneuos Neumann condition would result in the first term being zero. 
Further contributions of the second term are only vanishing for partially tangential $\ext{\qb}|_{\partial \manifold}$  
(\ie\ only the fully normal contribution $ (\normal\iP{}\ext{\qb}\iP{2}\normal)|_{\partial\manifold} $ is allowed to be nonzero), or 
the curvature is vanishing at $ \partial\manifold $ in direction of the conormal $ \conormal $, generally for all $ \ext{\psib} $ 
and $ \partial\manifold\neq\emptyset $.

\section{Identities for geometric Quantities}
\label{sec:geometricIdentities}
For simplifying the weak formulation of vector- and tensor-valued Laplacians we used the following identities
\begin{align}
0 & = \shapeOperator^2 - \meanCurvature \shapeOperator + \gaussianCurvature \Pi \\
\| \shapeOperator \|^2 & = \meanCurvature^2 - 2\gaussianCurvature \label{eq:shapesquare}\\
\shapeOperator \inP \Pi & = \shapeOperator.
\end{align}

\section{Thin film Christoffel symbols}
\label{sec:AppendixThinFilmChristoffelSymbols}
To define the metric compatible derivation also in the thin film we define thin film Christoffel symbols along the choice of coordinates. To clearly distinguish between the usage of thin film and manifold coordinates we use here upper case index letters to denote the thin film coordinates. The Christoffel symbols (of second kind) in the thin film are
  \begin{align}
    \ChristE{K}{IJ} &= \frac{1}{2} G^{KL}\left( \partial_{I}G_{JL} + \partial_{J}G_{IL} - \partial_{L}G_{IJ} \right)\formPeriod
  \end{align}
For the thin film metric \( \Gb \) mixed tangential-normal components are zero (which also holds for the inverse metric) and the pure normal component is constant. Hence, we obtain
    \begin{align}
      \ChristE{K}{\xi\xi} &= \frac{1}{2}G^{KL}\left( \partial_{\xi}G_{\xi L} + \partial_{\xi}G_{\xi L} - \partial_{L}G_{\xi\xi} \right)
                        = 0 \formComma\\
      \ChristE{\xi}{I\xi} &= \frac{1}{2} G^{\xi \xi}\left( \partial_{I}G_{\xi\xi} + \partial_{\xi}G_{I \xi} - \partial_{\xi}G_{I \xi} \right)
                        = 0\formPeriod
    \end{align}
The partial derivative in normal direction for the tangential part of the thin film metric is
	\begin{align}
	  \partial_{\xi}G_{ij} &= 2\left( -B_{ij} + \xi\left[ \shapeOperator^{2} \right]_{ij} \right)\formPeriod
	\end{align}
The pure tangential part of the inverse metric is $ G^{ij} = [(\gb - \xi\shapeOperator)^{-2}]^{ij} $ and can be approximated by Taylor or Neumann series, 
see \cite{Nitschkeetal_PRSA_2018}.
Therefore, we obtain
	\begin{align}
	  \ChristE{\xi}{ij} &= \frac{1}{2} G^{\xi\xi}\left( \partial_{i}G_{j\xi} + \partial_{j}G_{i\xi} - \partial_{\xi}G_{ij} \right)
	                 = B_{ij} - \xi\left[ \shapeOperator^{2} \right]_{ij} \formComma \\
	  \ChristE{k}{i\xi} &= \frac{1}{2} G^{kl}\left( \partial_{i}G_{\xi l} + \partial_{\xi}G_{il} - \partial_{l}G_{i\xi} \right)
	                 = \left( g^{kl} + 2\xi B^{kl} +  \landau(\xi^{2}) \right)\left( -B_{il}  + \xi\left[ \shapeOperator^{2} \right]_{il} \right)\\
	                &= -\tensor{B}{_{i}^{k}} - \xi\tensor{\left[ \shapeOperator^{2} \right]}{_{i}^{k}} + \landau(\xi^{2}) 
	\end{align}
 and with the substitution \eqref{eq:shapesquare} the remaining two terms.
 For the pure tangential thin film Christoffel symbols, we introduce the auxillary notion of $ \tensor{\beta}{_{ij}^{k}} := \nabla_j\tensor{B}{_{i}^{k}} +  \nabla_i\tensor{B}{_{j}^{k}} - \nabla^k\tensor{B}{_{ij}}$. Expressing \( \tensor{\beta}{_{ij}^{k}} \) in terms of partial derivatives and use the symmetry of the shape operator, we obtain
 \begin{align}
   \tensor{\beta}{_{ij}^{k}} 
                &= g^{kl}\left( B_{il|j} + B_{jl|i} - B_{ij|l} \right) \\
                &= g^{kl}\left( \partial_{j}B_{il} - \Christ{ij}{m}B_{ml} - \Christ{jl}{m}B_{im} 
                              + \partial_{i}B_{jl} - \Christ{ij}{m}B_{ml} - \Christ{il}{m}B_{jm} 
                              - \partial_{l}B_{ij} + \Christ{il}{m}B_{mj} + \Christ{jl}{m}B_{im}\right)\\
                &= g^{kl}\left( \partial_{j}B_{il} + \partial_{i}B_{jl} - \partial_{l}B_{ij} - 2\Christ{ij}{m}B_{ml}  \right) \\
                &= g^{kl}\left( \partial_{j}B_{il} + \partial_{i}B_{jl} - \partial_{l}B_{ij} \right)
                     -2\Gamma_{ijl}B^{kl} \\
                &= g^{kl}\left( \partial_{j}B_{il} + \partial_{i}B_{jl} - \partial_{l}B_{ij} \right)
                     -B^{kl}\left( \partial_{j}g_{il} + \partial_{i}g_{jl} - \partial_{l}g_{ij} \right) \formPeriod
 \end{align}
Using this expression we yield for the pure tangential parts of the thin film Christoffel symbols
 \begin{align}
   \ChristE{k}{ij} &= \frac{1}{2} G^{kl}\left( \partial_{i}G_{jl} + \partial_{j}G_{il} - \partial_{l}G_{ij} \right) \\
               &=  \frac{1}{2} \left( g^{kl} + 2\xi B^{kl} +  \landau(\xi^{2}) \right)
                     \left( \partial_{j}g_{il} + \partial_{i}g_{jl} - \partial_{l}g_{ij} 
                               -2\xi\left( \partial_{j}B_{il} + \partial_{i}B_{jl} - \partial_{l}B_{ij} \right)\right)\\
               &= \Christ{ij}{k} + \xi\left( B^{kl}\left( \partial_{j}g_{il} + \partial_{i}g_{jl} - \partial_{l}g_{ij} \right) 
                                            - g^{kl}\left( \partial_{j}B_{il} + \partial_{i}B_{jl} - \partial_{l}B_{ij} \right)\right)
                               +  \landau(\xi^{2})\\
               &= \Christ{ij}{k} - \xi\tensor{\beta}{_{ij}^{k}} + \landau(\xi^{2})\formPeriod
 \end{align}
Summarizing these results we can describe the thin film Christoffel symbols by a second order expansions in normal direction of the manifold Christoffel symbols:
    \begin{align}
      \ChristE{k}{ij} &= \Christ{ij}{k} - \xi\tensor{\beta}{_{ij}^{k}} + \landau(\xi^{2})\\
      \ChristE{\xi}{ij} &= B_{ij} - \xi \left[ \shapeOperator^{2} \right]_{ij}  &
      \ChristE{k}{i\xi} = \ChristE{k}{\xi i} &= -\tensor{B}{_{i}^{k}} - \xi\tensor{\left[ \shapeOperator^{2} \right]}{_{i}^{k}} + \landau(\xi^{2}) \\
      \ChristE{K}{\xi\xi} &=  \ChristE{\xi}{I\xi} = \ChristE{\xi}{\xi I} = 0 \formPeriod
    \end{align}

\section{First order differential operators on manifolds}
\label{sec:AppendixDifferentialOperators}
\paragraph{Divergence of tensor field}
For a tensor field $\tb \in \tensorField{d}{\manifold}$ we can define canonical a divergence by 
\begin{align}
\Div \tb = \nabla \tb \iP{d,d+1} \gb \in \tensorField{d-1}{\manifold}
\end{align}
Using the thin film extension we can describe $\gb_{ij}=\ProjSurf_{ij}$ and insert the Cartesian description of the gradient $\nabla \tb$ to yield
\begin{align}
[\Div \tb ]_{I_1 \hdots I_{d-1}}= & [\nabla\tb]_{I_1 \hdots I_d K} \tensor{\ProjSurf}{^{I_d K}} \nonumber \\
 = & \ProjSurf[\extNabla\ext{\tb}]_{I_1 \hdots I_d K}\tensor{\ProjSurf}{^{I_d K}} + [\tensor{B}{_{K I_1}} \tensor{\ext{t}}{^{L}_{I_2 \hdots I_d}} \nu_L +  \hdots + \tensor{B}{_{K I_d}} \tensor{\ext{t}}{_{I_1 \hdots I_{d-1}}^{L}} \nu_L ] \tensor{\ProjSurf}{^{I_d K}} \nonumber \\
 = & \tensor{\ProjSurf[\extNabla\ext{\tb}]}{_{I_1 \hdots I_{d-1}}^{K}_{K}} + [\tensor{B}{^{I_d}_{I_1}} \tensor{\ext{t}}{^{L}_{I_2 \hdots I_d}} \nu_L +  \hdots + \tensor{B}{^{I_d}_{I_d}} \tensor{\ext{t}}{_{I_1 \hdots I_{d-1}}^{L}} \nu_L ]
\end{align}

\paragraph{Rotation of tensor field}
The rotation of a tensor field $\tb \in \tensorField{d}{\manifold}$ can not be defined in a canonical way, but the definitions can be clustered into two groups. One type of definition increases the tensorial degree of $\tb$ while the second type decreases the tensorial degree. Anyhow both definitions include the contraction of $\nabla \tb$ with the Levi-Civita tensor $\boldsymbol{E}\in\tensorField{n}{\manifold}$, but vary in the amount of components and which components are contracted.
The Levi-Civita Tensor is given uniquely by its fully pairwise skew-symmetric behavior and normalization condition 
$ \left\| \boldsymbol{E} \right\|^2 = \left\| \gb \right\|^2 = n $ and can be obtained by the tonsorial description of the volume form $ \mathrm{d}\manifold $,
a differential $ n $-form, \ie\ $ E_{i_1\hdots i_n} = \mathrm{d}\manifold\left( \partial_{i_1}\Xc,\hdots, \partial_{i_n}\Xc \right) $.
We present here two examples for tensor fields on a $ n $-dimensional manifold
\begin{align}
\Rot \tb =  \nabla \tb \iP{d+1} \boldsymbol{E} \in \tensorField{d+n-1}{\manifold}  \quad \mbox{ and } \quad \rot_k \tb = -  \nabla \tb \iP{k, d+1} \boldsymbol{E} \in \tensorField{d+n - 3}{\manifold}\formPeriod
\end{align}
Using for example the Cartesian description of the covariant derivative and a thin film extension of $\boldsymbol{E}$ we yield 
\begin{align}
[\Rot \tb]_{I_1 \hdots I_d J_1 \hdots J_{n-1}}& = [\nabla \tb ]_{I_i\hdots I_d L} \tensor{E}{^{L}_{J_1 \hdots J_{n-1}}} \nonumber \\
& = \ProjSurf[\extNabla\ext{\tb}]_{I_1 \hdots I_d L}\tensor{E}{^{L}_{J_1 \hdots J_{n-1}}} + [\tensor{B}{_{L I_1}} \tensor{\ext{t}}{^{M}_{I_2 \hdots I_d}} +  \hdots + \tensor{B}{_{L I_d}} \tensor{\ext{t}}{_{I_1 \hdots I_{d-1}}^{M}} ]\nu_M \tensor{E}{^{L}_{J_1 \hdots J_{n-1}}}
\end{align}
and
\begin{align}
&[\rot_k \tb]_{I_1 \hdots I_{d-1}J_1 \hdots J_{n-2} } \nonumber \\
= & -[\nabla \tb ]_{I_1\hdots I_{k-1} K I_{k} \hdots I_{d-1} L} \tensor{E}{^{KL}_{J_1 \hdots J_{n-2}}}\notag \\
     =&  \ProjSurf[\extNabla\ext{\tb}]_{I_1\hdots I_{k-1} K I_{k} \hdots I_{d-1} L} \tensor{E}{^{LK}_{J_1 \hdots J_{n-2}}} 
    + [\tensor{B}{_{L I_1}} \tensor{\ext{t}}{^{M}_{I_2 \hdots I_{k-1} K I_{k} \hdots I_{d-1}}}    \\
    &\quad\quad\quad\quad +\hdots 
      +\tensor{B}{_{L K}} \tensor{\ext{t}}{_{I_1\hdots I_{k-1}}^{M}_{I_{k} \hdots I_{d-1}}} +
     \hdots + \tensor{B}{_{L I_{d-1}}} \tensor{\ext{t}}{_{I_1\hdots I_{k-1} K I_{k} \hdots I_{d-2}}^{M}} ]\nu_{M} \tensor{E}{^{LK}_{J_1 \hdots J_{n-2}}} \notag 
    \formPeriod
\end{align}

\bibliographystyle{siam}

\begin{thebibliography}{10}

\bibitem{abraham2012manifolds}
{\sc R.~Abraham, J.~E. Marsden, and T.~Ratiu}, {\em Manifolds, tensor analysis,
  and applications}, vol.~75, Springer Science \& Business Media, 2012.

\bibitem{Andrews2011}
{\sc B.~Andrews and C.~Hopper}, {\em The Ricci Flow in Riemannian Geometry: A
  Complete Proof of the Differentiable 1/4-Pinching Sphere Theorem}, no.~Nr.
  2011 in Lecture Notes in Mathematics, Springer, 2011.

\bibitem{Backus_ARMA_1966}
{\sc G.~E. Backus}, {\em Potentials for tangent tensor fields on spheroids},
  Arch. Ration. Mech. Anal., 22 (1966), pp.~210--252.

\bibitem{Barreraetal_EJP_1985}
{\sc R.~G. Barrera, G.~A. Estevez, and J.~Giraldo}, {\em Vector spherical
  harmonics and their application to magnetostatics}, Eur. J. Phys., 6 (1985),
  p.~287.

\bibitem{Bertalmioetal_JCP_2001}
{\sc M.~Bertalmio, L.~T. Cheng, S.~Osher, and G.~Sapiro}, {\em Variational
  problems and partial differential equations on implicit surfaces}, J. Comput.
  Phys., 174 (2001), pp.~759--780.

\bibitem{Burmanetal_IJNME_2015}
{\sc E.~Burman, S.~Claus, P.~Hansbo, M.~Larson, and A.~Massing}, {\em Cutfem:
  Discretizing geometry and partial differential equations}, Int. J. Numer.
  Meth. Engng., 104 (2015), pp.~472--501.

\bibitem{Craneetal_ACM_2013}
{\sc K.~Crane, F.~de~Goes, M.~Desbrun, and P.~Schr\"oder}, {\em {Digital
  geometry processing with discrete exterior calculus}}, ACM SIGGRAPH 2013, 7
  (2013).

\bibitem{Demlow_SIAMJNA_2009}
{\sc A.~Demlow}, {\em {Higher-order finite element methods and pointwise error
  estimates for elliptic problems on surfaces}}, SIAM J. Numer. Anal., 47
  (2009), pp.~805--827.

\bibitem{DziukElliott_IMAJNA_2007}
{\sc G.~Dziuk and C.~M. Elliott}, {\em Finite elements on evolving surfaces},
  IMA J. Num. Ana., 27 (2007), p.~261.

\bibitem{DziukElliott_JCM_2007}
\leavevmode\vrule height 2pt depth -1.6pt width 23pt, {\em Surface finite
  elements for parabolic equations}, J. Comput. Math., 25 (2007), p.~385.

\bibitem{Dziuketal_IFB_2008}
\leavevmode\vrule height 2pt depth -1.6pt width 23pt, {\em Eulerian finite
  element method for parabolic {PDE}s on implicit surfaces}, Interf. Free
  Bound., 10 (2008), p.~119.

\bibitem{Dziuk}
\leavevmode\vrule height 2pt depth -1.6pt width 23pt, {\em {Finite element
  methods for surface {PDE}s}}, Acta Numerica, 22 (2013), p.~289.

\bibitem{Freedenetal_MG_1995}
{\sc W.~Freeden, T.~Gervens, and M.~Schreiner}, {\em Tensor spherical harmonics
  and tensor spherical splines}, Manuscr. Geodaet., 19 (1994), pp.~80--100.

\bibitem{Freedenetal_Springer_2009}
{\sc W.~Freeden and M.~Schreiner}, {\em Spherical Functions of Mathematical
  Geosciences -- A Scalar, Vectorial, and Tensorial Setup}, Advances in
  Geophysical and Environmental Mechanics and Mathematics, Springer, 2009.

\bibitem{DeGoesetal_ACM_2015}
{\sc F.~D. Goes, M.~Desbrun, and Y.~Tang}, {\em {Vector field processing on
  triangle meshes}}, ACM SIGGRAPH Asia 2015, 17 (2015).

\bibitem{Greeretal_JCP_2006}
{\sc J.~Greer, A.~L. Bertozzi, and G.~Sapiro}, {\em Fourth order partial
  differential equations on general geometries}, J. Comput. Phys., 216 (2006),
  p.~216.

\bibitem{Grossetal_JCP_2018}
{\sc B.~J. Gross and P.~J. Atzberger}, {\em {Hydrodynamic flows on curved
  surfaces: Spectral numerical methods for radial manifold shapes}}, J. Comput.
  Phys., 371 (2018), pp.~663--689.

\bibitem{Grossetal_arXiv_2017}
{\sc S.~Gro{\ss}, T.~Jankuhn, M.~A. Olshanskii, and A.~Reusken}, {\em {A trace
  finite element method for vector-Laplacians on surfaces}}, arXiv:1709.00479.

\bibitem{Hansboetal_arXiv_2017}
{\sc P.~Hansbo, M.~G. Larson, and K.~Larsson}, {\em {Analysis of finite element
  methods for vector Laplacians on surfaces}}, arXiv:1610.06747.

\bibitem{Hirani_ACM_2003}
{\sc A.~N. Hirani}, {\em {Discrete exterior calculus}}, ACM Doctoral
  Dissertation,  (2003).

\bibitem{Jahnkuhnetal_arXiv_2017}
{\sc T.~Jankuhn, M.~A. Olshanskii, and A.~Reusken}, {\em {Incompressible fluid
  problems on embedded surfaces: Modeling and variational formulations}},
  arXiv:1702.02989.

\bibitem{Jaenich2013}
{\sc L.~Kay and K.~J{\"a}nich}, {\em Vector Analysis}, Undergraduate Texts in
  Mathematics, Springer New York, 2013.

\bibitem{Keberetal_Science_2014}
{\sc F.~C. Keber, E.~Loiseau, T.~Sanchez, S.~J. DeCamp, L.~Giomi, M.~Bowick,
  M.~C. Marchetti, Z.~Dogic, and A.~R. Bausch}, {\em Topology and dynamics of
  active nematic vesicles}, Science, 345 (2014), p.~1135.

\bibitem{kuhnel2006differential}
{\sc W.~K{\"u}hnel, B.~Hunt, and A.~M. Society}, {\em Differential Geometry:
  Curves - Surfaces - Manifolds}, Student mathematical library, American
  Mathematical Society, 2006.

\bibitem{Dunkeletal_PRL_2018}
{\sc O.~Mickelin, J.~Slomka, K.~J. Burns, D.~Lecoanet, G.~M. Vasil, L.~M.
  Faria, and J.~Dunkel}, {\em {Anomalous chained turbulence in actively driven
  flows on spheres}}, Phys. Rev. Lett., 120 (2018), p.~164503.

\bibitem{Mohamedetal_JCP_2016}
{\sc M.~S. Mohamed, A.~N. Hirani, and R.~Samtaney}, {\em {Discrete exterior
  calculus discretization of incompressible Navier-Stokes equations over
  surface simplicial meshes}}, J. Comput. Phys., 312 (2016), p.~175.

\bibitem{Nestleretal_JNS_2018}
{\sc M.~Nestler, I.~Nitschke, S.~Praetorius, and A.~Voigt}, {\em {Orientational
  order on surfaces - the coupling of topology, geometry and dynamics}}, J.
  Nonlin. Sci., 28 (2018), p.~147.

\bibitem{Nestleretal_arXiv_2018}
{\sc M.~Nestler, S.~Reuther, and A.~Voigt}, {\em {Hydrodynamic interactions in
  polar liquid crystals on evolving surfaces}}, arXiv,  (2018).

\bibitem{Nitschkeetal_PRSA_2018}
{\sc I.~Nitschke, M.~Nestler, S.~Praetorius, H.~L\"oowen, and A.~Voigt}, {\em
  {Nematic liquid crystals on curved surfaces - a thin film limit}}, Proc.
  Royal Soc. London A, 474 (2018), p.~2214.

\bibitem{Nitschkeetal_book_2017}
{\sc I.~Nitschke, S.~Reuther, and A.~Voigt}, {\em {Discrete Exterior Calculus
  (EC) for the Surface Navier-Stokes Equation}}, in Transport Processes at
  Fluidic Interfaces, D.~Bothe and A.~Reusken, eds., Springer, 2017,
  pp.~177--197.

\bibitem{Nitschkeetal_arXiv_2018}
{\sc I.~Nitschke and A.~Voigt}, {\em {in preperation}}.

\bibitem{Nitschkeetal_JFM_2012}
{\sc I.~Nitschke, A.~Voigt, and J.~Wensch}, {\em {A finite element approach to
  incompressible two-phase flow on manifolds}}, J. Fluid Mech., 708 (2012),
  p.~418.

\bibitem{Olshanskiietal_arXiv_2018}
{\sc M.~Olshanskii, A.~Quaini, A.~Reusken, and V.~Yushutin}, {\em {A finite
  element method for the surface Stokes problem}}, arXiv:1801.06589.

\bibitem{Olshanskiietal_LNCSE_2018}
{\sc M.~Olshanskii and A.~Reusken}, {\em {Trace finite element methods for PDEs
  on surfaces}}, in Geometrically unfitted finite element methods and
  applications, S.~Bordas, E.~Burman, M.~Larson, and M.~Olshanskii, eds.,
  Springer, 2018, pp.~211--258.

\bibitem{Praetoriusetal_PRE_2018}
{\sc S.~Praetorius, A.~Voigt, R.~Wittkowski, and H.~L\"owen}, {\em {Active
  crystals on a sphere}}, Phys. Rev. E, 97 (2018), p.~052615.

\bibitem{Raetzetal_CMS_2006}
{\sc A.~R{\"a}tz and A.~Voigt}, {\em {PDE}'s on surfaces---a diffuse interface
  approach}, Commun. Math. Sci., 4 (2006), pp.~575--590.

\bibitem{Reutheretal_PF_2018}
{\sc S.~Reuther and A.~Voigt}, {\em {Solving the incompressible surface
  Navier-Stokes equation by surface finite elements}}, Phys. Fluids, 30 (2018),
  p.~012107.

\bibitem{Schouten1954}
{\sc J.~Schouten}, {\em Ricci-calculus, an introduction to tensor analysis and
  its geometrical applications, by J.A. Schouten, ... 2nd edition ...},
  Springer-Verlag, 1954.

\bibitem{Stoeckeretal_JIS_2008}
{\sc C.~St{\"o}cker and A.~Voigt}, {\em Geodesic evolution laws - a level set
  approach}, SIAM J. Imag. Sci., 1 (2008), p.~379.

\bibitem{Unal1995}
{\sc B.~{\"U}nal}, {\em Divergence theorems in semi-riemannian geometry}, Acta
  Applicandae Mathematica, 40 (1995), pp.~173--178.

\bibitem{Veyetal_CVS_2007}
{\sc S.~Vey and A.~Voigt}, {\em {AMDiS}: adaptive multidimensional
  simulations}, Comput. Vis. Sci., 10 (2007), pp.~57--67.

\bibitem{Witkowskietal_ACM_2015}
{\sc T.~Witkowski, S.~Ling, S.~Praetorius, and A.~Voigt}, {\em Software
  concepts and numerical algorithms for a scalable adaptive parallel finite
  element method}, Adv. Comput. Math., 41 (2015), pp.~1145--1177.

\end{thebibliography}

\end{document}